\documentclass[twoside]{amsbook}
\usepackage[french,german,english]{babel} 
\usepackage[latin1]{inputenc} 
\usepackage{amsmath}
\usepackage{amsfonts} 
\usepackage{amsthm}
\usepackage{amscd}
\usepackage{amsopn}
\usepackage{amstext}
\usepackage{amsxtra}
\usepackage{textcase} 
\usepackage[ruled,linesnumbered,algochapter]{algorithm2e} 
\usepackage{bm} 
\numberwithin{equation}{section}
\numberwithin{figure}{chapter}

\newcounter{algorithm}[chapter]

\theoremstyle{plain}
\newtheorem{theorem}{Theorem}[chapter]

\newtheorem{cl}[theorem]{Corollary}
\newtheorem{pr}[theorem]{Proposition}

\theoremstyle{definition}
\newtheorem{deff}[theorem]{Definition}

\theoremstyle{remark}
\newtheorem{note}[theorem]{Note}
\newtheorem{remark}[theorem]{Remark}

\providecommand{\abs}[1]{\left\lvert #1 \right\rvert}

\providecommand{\set}[1]{\left\lbrace #1 \right\rbrace}
\providecommand{\gen}[1]{\left\langle #1 \right\rangle}
\providecommand{\Sym}[1]{\operatorname{Sym}( #1 )}
\renewcommand{\Pr}[1]{\operatorname{Pr}[ #1 ]}

\newcommand{\field}[1]{\mathbb{#1}}

\newcommand{\F}{\field{F}}
\newcommand{\GAP}{\textsf{GAP}}

\newcommand{\C}{\mathcal{C}}

\DeclareMathOperator{\GL}{GL}
\DeclareMathOperator{\GF}{GF}

\DeclareMathOperator{\SL}{SL}
\DeclareMathOperator{\GO}{GO}
\DeclareMathOperator{\SO}{SO}
\DeclareMathOperator{\Sz}{Sz}
\DeclareMathOperator{\Ree}{Ree}

\newcommand{\OR}[1]{\operatorname{O} ( #1 )}

\begin{document}

\frontmatter

\title{The Schreier-Sims algorithm for matrix groups}
\author{Henrik B\"a\"arnhielm}
\date{2004-07-17}


\begin{abstract}
  This is the report of a project with the aim to make a new
  implementation of the Schreier-Sims algorithm in \GAP, specialized
  for matrix groups. The standard Schreier-Sims algorithm is
  described in some detail, followed by descriptions of the probabilistic
  Schreier-Sims algorithm and the Schreier-Todd-Coxeter-Sims
  algorithm. Then we discuss our implementation and some
  optimisations, and finally we report on the performance of our
  implementation, as compared to the existing implementation in \GAP,
  and we give benchmark results. The conclusion is that our
  implementation in some cases is faster and consumes much less memory.
\end{abstract}

\urladdr{http://matrixss.sourceforge.net/}
\email{henrik.baarnhielm@imperial.ac.uk}

\maketitle



\tableofcontents


\chapter*{Preface}

This report is made as part of a student project at the Department of
Mathematics, Imperial College of Science, Technology and Medicine in
London, United Kingdom. It is made in partial fulfillment of the
requirements for the degree of Master of Science in Pure Mathematics.

I would like to thank my supervisors Prof. Alexander
Ivanov at Imperial College and Dr. Leonard Soicher at Queen Mary,
University of London, for their encouragement and help.

I am also particularly grateful to Alexander Hulpke at Colorado
State University in Fort Collins, USA, for suggesting this project to
me, and for his constant help and advice during my work.

\mainmatter

\chapter{Introduction}
The following report makes up one the two parts of a project in
computational group theory, the other part being a software package
for the computer system \GAP~(see \cite{GAP}), which can be found on
the WWW at the URL given on the title page. This text will describe
the mathematics that provide the foundations of the package, including
the algorithms used and their complexity, and also some of the more
computer science oriented aspects, like what datastructures that were
used, and how the implementation was done.

Computational group theory (CGT) is an area of research on the border
between group theory and computer science, and work in CGT is often of
both theoretical (mathematical) and practical (programming) nature,
leading to both theoretical results (mathematical theorems and proofs)
and practical results (software), and this project is no exception.
Introductory surveys of CGT can be found in \cite{sims98},
\cite{seress97}, \cite{neubuser95} and \cite{cannon92}.

The aim of the project was to make a \GAP~package with an
implementation of the Schreier-Sims algorithm for matrix groups. The
Schreier-Sims algorithm computes a base and a strong generating set
for a group, and an implementation of this fundamental algorithm is
already included in the standard \GAP~distribution, but that
implementation always first computes a faithful action (ie a
permutation representation) of the given group and then executes the
algorithm on the resulting permutation group. The idea for this
project was to restrict attention to matrix groups, and implement a
version of the algorithm which works with the matrices directly, and
see if one can obtain a more efficient implementation in this way.

A survey of computational matrix group theory can be found in
\cite{niemeyer01}. It should be noted that we are only interested in
finite groups, ie matrix groups over finite fields, and therefore
there is no need to worry about any noncomputability or undecidability
issues.

We will begin with a quick reference of the basic concepts from group
theory and computer science that are being used, before moving on to
describe the Schreier-Sims algorithm. The description will be quite
detailed, and then we will turn to the variants of the algorithm that
have also been implemented in the project: the random (ie.
probabilistic) Schreier-Sims algorithm and the
Schreier-Todd-Coxeter-Sims algorithm. After that we will say something
about the implementation, and describe some tricks and
improvements that have been done to make the algorithm faster.
Finally, the practical performance and benchmark results of the
implementation will be shown, and compared to the existing
implementation in \GAP.

It must be mentioned that a report similar to this one is
\cite{murray93}, from which a fair amount of inspiration comes.

\chapter{Preliminaries}
The definitions and statements in this section are assumed to be
known, but we state them anyway, since authors often use different
notation and sometimes put slightly different meaning to some of the
following concepts (eg. the exact definition of graphs tend to
vary). 

\section{Graph theory}
First some graph theory, where we follow \cite{biggs89}.

\begin{deff} \label{def_graph}
A \emph{directed graph} is an ordered pair $G = (V, E)$ where $V$ is a finite non-empty set, the \emph{vertices} of $G$ and $E \subseteq V \times V$ is the \emph{edges} of $G$.
\end{deff}
\begin{remark}
  A graph in the sense of \ref{def_graph} is sometimes called a
  \emph{combinatorial graph} in the literature, to emphasize that they
  are not \emph{metric graphs} in the sense of \cite{bridson99}. As we
  are interested only in finite graphs and not in geometry, we do not
  make use of this nomenclature.
\end{remark}

\begin{remark}
The definition implies that our graphs have no multiple edges, but may have loops. We will henceforth omit the word ''directed'' since these are the only
graphs we are interested in.
\end{remark}

\begin{deff} \label{def_graph_misc}
If $G = (V, E)$ is a graph, a sequence $v = v_1, v_2, \dotsc, v_n =
  u$ of vertices of $G$ such that $(v_i, v_{i + 1}) \in E$ for $i = 1,
  \dotsc, n - 1$ is called a \emph{walk} from $v$ to $u$. The \emph{length} of a walk $v = v_1, \dotsc, v_k = u$ is $k - 1$. If $v = u$
  then the walk is called a \emph{cycle}. If all vertices in the walk are
  distinct, it is called a \emph{path}. If every pair of distinct vertices in $G$ can be joined by a path, then $G$ is \emph{connected}.
\end{deff}

\begin{deff}
A graph $G^{\prime} = (V^{\prime}, E^{\prime})$ is a \emph{subgraph} of a graph $G = (V, E)$ if $V^{\prime} \subseteq V$ and $E^{\prime} \subseteq E$.
\end{deff}

\begin{deff} \label{def_tree}
A \emph{tree} is a connected graph without cycles. A \emph{rooted tree} is a tree $T = (V, E)$ where a vertex $r \in V$ have been specified as the \emph{root}.
\end{deff}

\begin{deff}
If $G = (V, E)$ is a graph, then a tree $T$ is a \emph{spanning tree} of $G$ if $T$ is a subgraph of $G$ and $T$ have vertex set $V$.
\end{deff}

The following are elementary and the proofs are omitted.
\begin{pr}
In a tree $T = (V, E)$ we have $\abs{E} = \abs{V} - 1$ and there
exists a \emph{unique} path between every pair of distinct
vertices. Conversely, if $G$ is a graph where every pair of distinct
vertices can be joined by a unique path, then $G$ is a tree.
\end{pr}

\begin{pr}
Every graph has a spanning tree.
\end{pr}

\begin{deff} 
Let $T = (V, E)$ be a rooted tree with root $r \in V$. The \emph{depth} of a node $x \in V$ is the length of the path from $x$ to $r$. The \emph{height} of $T$ is the maximum depth of any node.
\end{deff}

\begin{deff}
If $G = (V, E)$ is a graph, then a \emph{labelling} of $G$ is a function $w : E \to L$ where $L$ is some set of ''labels''. A \emph{labelled graph} is a graph with a corresponding labelling.
\end{deff}

\section{Group theory}
Now some group theory, where we follow our standard references \cite{bb96} and \cite{rose78}.
\begin{note}
All groups in this report are assumed to be finite.
\end{note}

\begin{deff} \label{def_cayley}
If $G = \gen{S}$ is a group, then the \emph{Cayley graph} $\C_G(S)$ is the graph with vertex set $G$ and edges $E = \set{ (g, sg) \mid g \in G, s \in S}$.
\end{deff}

\begin{deff}
An \emph{action} of a group $G$ on a finite set $X$ is a homomorphism
$\lambda : G \to \Sym{X}$ (where $\Sym{X}$ is the group of
permutations on $X$). If $\lambda$ is injective, the action is
\emph{faithful}. 
\end{deff}
\begin{remark}
Following a convention in computational group theory, actions are from the right, and $\lambda(g)x$ is abbreviated with $x^g$, for $g \in G$ and $x \in X$. The elements of $X$ are called \emph{points}. Note that some of the rules for exponents hold since we have an action, eg $(x^g)^h = x^{gh}$.
\end{remark}

\begin{deff}
Let $G$ be a group acting on the finite set $X$. For each point
$\alpha \in X$, the \emph{orbit} of $\alpha$ is $\alpha^G = \set{\beta
\in X \mid \beta = \alpha^g, g \in G}$ and the \emph{stabiliser} of
$\alpha$ is $G_{\alpha} = \set{g \in G \mid \alpha^g = \alpha}$.
\end{deff}

The following are elementary and the proofs are omitted.
\begin{pr}
Let $G$ be a group acting on the finite set $X$. For each $p \in X$ we have $G_p \leq G$, and so we can define inductively 
\begin{equation}
G_{\alpha_1, \alpha_2, \dotsc, \alpha_n} = (G_{\alpha_1, \alpha_2, \dotsc, \alpha_{n - 1}})_{\alpha_n}
\end{equation}
where $n > 1$ and $\alpha_1, \dotsc, \alpha_n \in X$.
\end{pr}

\begin{pr} \label{thm_orbit_stab}
Let $G$ be a group acting on the finite set $X$. For each $p \in X$, the map $\mu_p : G / G_p \to p^G$ given by
\begin{equation}
G_p g \mapsto p^g
\end{equation}
for each $g \in G$, is a bijection. In particular, $\abs{p^G} = [G : G_p]$.
\end{pr}

\section{Computer science}
When it comes to computer science, our standard reference is
\cite{clr90} where all basic computer science notions that are used here can be
found. First of all, complexity analysis of algorithms and its
asymptotic notation, in particular the $\OR{\cdot}$-notation, is assumed to be
known. Basic graph algorithms like breadth-first search, computation
of connected compontents and spanning tree algorithms are also assumed to be known.

Hash tables will likewise be used without further explanation. Even
though the code in the project does not use hashing explicitly, but
rely on \GAP~for that, it is worth mentioning that, to the author's
knowledge, the best general-purpose hash function known to humanity is
described in \cite{jenkins97}.

\chapter{The Schreier-Sims algorithm}

We will now describe the Schreier-Sims algorithm, and the references we are using are \cite{soicher98}, \cite{seress03} and \cite{butler91}. First we have to describe some background and clarify the problem that the algorithm solves.

\section{Background and motivation} \label{ss_background}
The overall goal for work in group theory is to understand groups, and
to answer various questions about groups. In particular, in CGT the
interest is focused on questions of algorithmic nature, so given a
group $G$ we are interested in algorithms for things like the following:
\begin{itemize}
\item What is $\abs{G}$?
\item List all elements of $G$, without repetition.
\item If $G \leq H$ and we are given an arbitrary $g \in H$, is it true that $g \in G$? This is referred to as \emph{the membership problem}.
\item Generate a random $g \in G$.
\item Is $G$ abelian (soluble, polycyclic, nilpotent)?
\item Find representatives for the conjugacy classes of $G$.
\item Find a composition series for $G$, including the (isomorphism classes of the) composition factors.
\item Given $g \in G$ or $H \leq G$ find the centralizer of $g$ or the normalizer of $H$, respectively.
\end{itemize}
To accomplish these tasks we need a computer representation of $G$, ie
a datastructure. A common way in which a group is given is via a
generating set, but this alone does not help in solving our problems,
so we need a better representation. It must be possible to compute
this representation from a generating set, and a nice property would
be if subgroups of $G$ in some direct way inherit the representation,
so that divide-and-conquer techniques could be used when designing
algorithms (see \cite{clr90}).

\section{Base and strong generating set} \label{section:bsgs}
Consider the situation where we have a chain of subgroups of $G$
\begin{equation}
G = G^0 \geq G^1 \geq \dotsb \geq G^n = 1
\end{equation}
Each $g \in G$ can be written $g = g_1 u_1$ where $u_1$ is a
representative of $G^1 g$ and $g_1 \in G^1$, and inductively we can
factorize $g$ as $g = u_n u_{n - 1} \dotsm u_1$, since the subgroup chain
reaches $1$. Moreover, this factorization is \emph{unique} since for each $i$ the cosets of $G^{i + 1}$ partition $G^i$, and by the same reason, different group elements will have different factorizations.

We thus see that if we know generating sets for the subgroups in such
a subgroup chain, and if we know a right transversal $T^i$ of the cosets of
$G^{i + 1}$ in $G^i$ for each $i = 0, \dotsc, n - 1$, then we could easily
solve at least the first two listed problems listed in section
\ref{ss_background}. By Lagrange we know $\abs{G} = \abs{T^1}
\abs{G^1}$ and inductively $\abs{G} = \abs{T^1} \dotsm \abs{T^n}$,
which solves the first problem. Using the factorization we have a
bijection from $G$ to $T^1 \times T^2 \times \dotsb \times T^n$ so by
enumerating elements of the latter set and multiplying in $G$ we can list the elements of $G$ without repetition.

Now we introduce a special type of subgroup chain.
\begin{deff}
  Let $G$ be a group acting on the finite set $X$. A sequence of points
  $(\alpha_1, \dotsc, \alpha_n)$ of $X$ such that
  $G_{\alpha_1, \dotsc, \alpha_n} = 1$ is called a \emph{base} for
  $G$. Note that the base determines a \emph{stabiliser chain}
\begin{equation}
G \geq G_{\alpha_1} \geq \dotsb \geq G_{\alpha_1, \dotsc, \alpha_n} = 1
\end{equation}
Let $G^i = G_{\alpha_1, \dotsc, \alpha_i}$ for $i = 1, \dotsc, n$. A generating set $S$ for $G$ such that $\gen{S \cap G^i} = G^i$ for all $i$ is called a \emph{strong generating set} (SGS) for $G$.
\end{deff}

The concept of a base and strong generating was first introduced in
\cite{sims70} in the context of permutation groups, and is of
fundamental importance, though mainly for permutation groups. We
already know that the first two of our problems can be solved if we
know a subgroup chain and we shall see that with a base and strong
generating set the membership problem can also easily be solved. The
Schreier-Sims algorithm computes a base and strong generating set for
a group $G$ given a generating set $S$, and since it is an efficient
algorithm if $G$ is a permutation group, the concept of base and
strong generating set has become very important. Many more
sophisticated algorithms for permutation groups require a base and
strong generating set as input. For matrix groups, on the other hand,
the situation is a bit more complicated, as we shall see later.

\section{Schreier trees}
Before giving the Schreier-Sims algorithm itself, there are a few
auxiliary algorithms that must be explained. Consider therefore a
group $G = \gen{S}$ that acts on a finite set $X$. From the previous
section we know that even if we have a base $(\alpha_1, \dotsc,
\alpha_n)$ for $G$, with corresponding stabiliser chain $G \geq G^1
\geq \dotsb G^n = 1$, we also need to find right transversals of the cosets of
$G^{i + 1}$ in $G^i$ for $i = 1, \dotsc, n - 1$. But since these groups
are stabilisers from the action on $X$, we can use Proposition
\ref{thm_orbit_stab} and instead find the orbits $\alpha_1^G,
\alpha_2^{G^1}, \dotsc, \alpha_n^{G^{n - 1}}$. We will see that the orbits are straightforward to compute.

The action of $G$ on $X$ can be represented by a labelled graph with
labels from $S$, analogous to the Cayley graph $\C_S(G)$ of $G$ from
Definition \ref{def_cayley}. Let the vertices of the graph be $X$ and the edges
be $\{(p, p^g) \mid p \in X, g \in G\}$, where the edge $(p, p^g)$ is
labelled by $g$. Obviously the orbits of the action are the connected
components of this graph, so $\alpha_1^G$ is the component containing
$\alpha_1$. We are interested in finding this orbit and to store it in
a datatstructure, and we are thus led to consider a spanning tree of
the connected component, since such a tree contains enough information
for us. The edges of the graph that are left out do not give us any
additional relevant information about the action of $G$, only
alternative ways to move between the points, and trees are
considerably easier to store than graphs.

\begin{deff} \label{def_bsgs}
  Let $G$ be a group acting on the finite set $X$ and let $\alpha \in
  X$. A spanning tree rooted at $\alpha$ for the component in the corresponding graph
  containing $\alpha$ is called a
  \emph{Schreier tree} for the orbit $\alpha^G$.
\end{deff}

The Schreier tree can be computed by a simple breadth-first search of
the component containing $\alpha$, and thus we have an algorithm for
finding the orbits. However, it is of course not a computationally
good idea to explictly generate the graph, then compute the connected
components and finally find the Schreier trees. A breadth-first search
to find the Schreier tree can be done without the graph itself, as
demonstrated by \ref{alg:orbit}.

\begin{algorithm} 
\dontprintsemicolon
\caption{\texttt{ComputeSchreierTree}}
\SetKwFunction{Tree}{Tree}
\SetKwFunction{AddChild}{AddChild}
\SetKwData{Points}{points}
\SetKwData{STree}{tree}
\SetKwData{Children}{children}
\KwData{A group $G = \gen{S}$ acting on a finite set $X$ and a point $\alpha \in X$.}
\KwResult{A Schreier tree for $\alpha^G$.}
\tcc{Assumes the existence of a function \texttt{Tree}$(x)$ that creates an empty tree with root $x$ and a function \texttt{AddChild}($T$, $p_1$, $p_2$, $l$) that adds $p_2$ as a child to $p_1$ in the tree $T$, with label $l$.}
\Begin{
  $\Points := \set{\alpha}$ \;
  $\STree := \Tree(\alpha)$ \;
  \Repeat{$\Points = \emptyset$}
  {
    $\Children := \emptyset$ \;
    \ForEach{$p \in \Points$}
    {
      \ForEach{$s \in S$}
      {
        $p^{\prime} := p^s$ \; \label{alg:orbit_innerloop1}
        \If{$p^{\prime} \notin \STree$}
        {
          $\AddChild(\STree, p, p^{\prime}, s)$ \;
          $\Children := \Children \cup \set{p^{\prime}}$ \; \label{alg:orbit_innerloop2}
        } 
      }
    }
    $\Points := \Children$ \;
  }
  \Return{$\STree$}
}
\refstepcounter{algorithm}
\label{alg:orbit}
\end{algorithm}
We noted above that we used Proposition \ref{thm_orbit_stab} to store the orbits
instead of the coset representatives, but we still need the latter, so
we must be able to compute them. Fortunately, that is straightforward to do when we
have a Schreier tree.  Assume we have a Schreier tree $T =
(V, E)$ for the orbit $\alpha^G$ for some point $\alpha \in X$. If $g
\in G$ then $\alpha^g \in V$ and there is a path from $\alpha$ to
$\alpha^g$ in $T$. Let $s_1, s_2, \dotsc, s_n$ be the labels of the
path, so that $s_i \in S$ for $i = 1, \dotsc, n$, and let $h = s_1 s_2
\dotsm s_n$. Then obviously $\alpha^h = \alpha^g$ so $G_{\alpha} g = 
G_{\alpha} h$ and $h$ is a coset representative for the coset of $g$.
Moreover, $h$ is unique since $T$ is a tree. Thus, to find the coset
representative for $g$ we only have to follow the unique path from
$\alpha^g$ to the root $\alpha$ and multiply the edge labels. \ref{alg:trace} performs the slightly more general task of
following the path from a given point to the root.

\begin{algorithm} 
\dontprintsemicolon
\caption{\texttt{OrbitElement}}
\SetKwFunction{EdgeLabel}{EdgeLabel}
\SetKw{KwRet}{break}
\SetKwData{Edge}{edge}
\KwData{A group $G = \gen{S}$ acting on a finite set $X$, a Schreier tree $T$ for the orbit $\alpha^G$ of the point $\alpha \in X$ and an arbitrary point $p \in X$.}
\KwResult{The element $g \in G$ such that $\alpha^g = p$}
\tcc{Assumes the existence of a function \texttt{EdgeLabel}$(T, p)$ that returns the label of the unique edge between $p$ and its parent in $T$}
\Begin
{
  $g := 1$ \;
  \While{$p \neq \alpha$}
  {
    $s := \EdgeLabel(T, p)$ \;
    $p := p^{s^{-1}}$ \;
    $g := s g$ \;
  }
  \Return{$g$}
}
\refstepcounter{algorithm}
\label{alg:trace}
\end{algorithm}

We defer the complexity analysis of these algorithms until later.

\section{Formulating the problem}
To more formally state the problem solved by the Schreier-Sims algorithm, the following is needed.
\begin{deff}
  Let $G$ be a group acting on the finite set $X$. A sequence of points
  $B = (\alpha_1, \dotsc, \alpha_n)$ of $X$ and a generating set $S$ for $G$, such that no element of $S$ fixes all points of $B$, is called a \emph{partial base} and \emph{partial strong generating set}, respectively. 
\end{deff}
\begin{remark}
A base and strong generating set as in Definition \ref{def_bsgs} are called \emph{complete}.
\end{remark}
\begin{remark} \label{rmk_partial_bsgs}
  If we define $G^i = G_{\alpha_1, \dotsc, \alpha_i}$, $S^i = S \cap
  G^i$ and $H^i = \gen{S^i}$ for $i = 1, \dotsc, n$, we see that $H^n
  = 1$ since no element of $S$ fixes all points of $B$ (and we use the
  convention that $\gen{\emptyset} = 1$). We therefore have
\begin{align}
G & \geq G^1 \geq \dotsb \geq G^n \\
G & \geq H^1 \geq \dotsb \geq H^n = 1 
\end{align}
Moreover, $G^i \geq \gen{S^i} = H^i$ for $i = 1, \dotsc, n$ and if we
have equality then, by Definition \ref{def_bsgs}, $S$ and $B$ are complete. If $h
\in H^{i + 1} = \gen{S^{i + 1}} = \gen{S \cap G^{i + 1}}$ then $h = s_1 \dotsm s_k$ where $\alpha_{i + 1}^{s_j} = \alpha_{i + 1}$ for $j = 1, \dotsc, k$ so $\alpha_{i + 1}^h = \alpha_{i + 1}$ and therefore $h \in H^i_{\alpha_{i + 1}}$. Thus $H^{i + 1} \leq H^i_{\alpha_{i + 1}}$ for $i = 0, \dotsc, n - 1$.
\end{remark}

Now our problem can be stated as follows: given a group $G$ acting on the
finite set $X$, together with a partial base $B$ with points from $X$
and partial strong generating set $S$, either verify that $B$ is a
(complete) base and that $S$ is a (complete) strong generating set, or
extend $B$ and $S$ so that they become complete. This is the problem that is
solved by the Schreier-Sims algorithm.

 The following result from \cite{leon80} is used when designing the algorithm.
\begin{theorem} \label{thm_leon}
  Let $G$ be a group acting on the finite set $X$, and let $B =
  (\alpha_1, \dotsc, \alpha_n)$ be a partial base and $S$ a partial
  strong generating set for $G$. Let also $G^i = G_{\alpha_1, \dotsc,
    \alpha_i}$, $S^i = S \cap G^i$, $H^i = \gen{S^i}$ for $i = 1,
  \dotsc, n$ and $G = G^0 = H^0$. Then the following statements are equivalent:
\begin{enumerate}
\item $B$ and $S$ are complete. \label{leon80_1}
\item $G^i = H^i$ for $i = 0, \dotsc, n$. \label{leon80_2}
\item $H^i_{\alpha_{i + 1}} = H^{i + 1}$ for $i = 0, \dotsc, n - 1$. \label{leon80_3}
\item $[H^i : H^{i + 1}] = \abs{\alpha_{i+1}^{H^i}}$ for $i = 0, \dotsc, n - 1$. \label{leon80_4}
\end{enumerate}
\end{theorem}
\begin{proof}
From Remark \ref{rmk_partial_bsgs} we know that \eqref{leon80_1} and \eqref{leon80_2} are equivalent. Assuming \eqref{leon80_2} we have
\begin{equation}
H^i_{\alpha_{i + 1}} = G^i_{\alpha_{i + 1}} = G^{i + 1} = H^{i + 1}
\end{equation}
for $i = 0, \dotsc, n - 1$ which is precisely \eqref{leon80_3}. If we instead assume \eqref{leon80_3} and also assume for induction that $G^i = H^i$ (the base case $G = H^0 = G^0$ is ok) then
\begin{equation}
G^{i + 1} = G^i_{\alpha_{i + 1}} = H^i_{\alpha_{i + 1}} = H^{i + 1}
\end{equation}
so by induction we get $G^i = H^i$ for $i = 0, \dotsc, n$, which is \eqref{leon80_2}.

Now assume \eqref{leon80_3} and note that from \ref{thm_orbit_stab} we
have $[H^i : H^i_{\alpha_{i + 1}}] = \abs{\alpha_{i+1}^{H^i}}$, so
since $H^i_{\alpha_{i + 1}} = H^{i + 1}$ we get \eqref{leon80_4}.
Finally, assume \eqref{leon80_4}. From Remark \ref{rmk_partial_bsgs} we know
$H^i_{\alpha_{i + 1}} \geq H^{i + 1}$ so if we again use
\ref{thm_orbit_stab} we get $\abs{\alpha_{i+1}^{H^i}} = [H^i : H^{i +
  1}] \geq [H^i : H^i_{\alpha_{i + 1}}] = \abs{\alpha_{i+1}^{H^i}}$
and thus $H^i_{\alpha_{i + 1}} = H^{i + 1}$.
\end{proof}

As observed earlier, we are often given a group in the form of a
generating set, but Schreier-Sims algorithm requires a partial base
and a partial strong generating set as input. Those are easy to
compute, though, using \ref{alg:partial_bsgs}. The algorithm also
makes sure that the partial strong generating set is closed under
inverses and does not contain the identity, which removes the need to
consider some special cases later on. We will see how the function
\texttt{NewBasePoint} that is used in \ref{alg:partial_bsgs} can be
implemented when $G$ is a matrix group.

\begin{algorithm} 
\dontprintsemicolon
\caption{\texttt{GetPartialBSGS}}
\SetKwFunction{NewBasePoint}{NewBasePoint}
\SetKwData{Base}{base}
\SetKwData{SGS}{sgs}
\SetKwData{Point}{point}
\KwData{A group $G = \gen{S}$ acting on a finite set $X$ and a sequence of points $B$ of $X$ (possibly empty).}
\KwResult{A partial base $B^{\prime}$ and partial strong generating set $S^{\prime}$ for $G$.}
\tcc{Assumes the existence of a function \texttt{NewBasePoint}$(g)$ that returns a point $p \in X$ such that $p^g \neq p$}
\Begin
{
  $\Base := B$ \;
  $\SGS := \emptyset$ \;
  \ForEach{$s \in S \setminus \set{1}$}
  {
    \If{$\Base^s = \Base$} 
    { \label{alg:partial_bsgs_innerloop1}
      $\Point := \NewBasePoint(s)$ \;
      $\Base := \Base \cup \set{\Point}$ \; \label{alg:partial_bsgs_innerloop2} 
    }
    $\SGS := \SGS \cup \set{s, s^{-1}}$ \;
  }    
  \Return{$(\Base, \SGS)$}
}
\refstepcounter{algorithm}
\label{alg:partial_bsgs}
\end{algorithm}

\section{Schreier's Lemma}
The name Schreier in Schreier-Sims algorithm comes from the following
result, which in our case allows us to find a generating set for a stabiliser. It
first appeared in \cite{schreier27}, and our proof is originally from \cite{hall59}.
\begin{theorem}[Schreier's Lemma] \label{schreiers_lemma}
Let $G = \gen{S}$ be a group, let $H \leq G$ and let $T$ be a right transversal of the cosets of $H$ in $G$. For $g \in G$, let $\bar{g} \in T$ be the unique element such that $Hg = H\bar{g}$. Then $H$ is generated by
\begin{equation}
S_H = \set{t s (\overline{ts})^{-1} \mid t \in T, s \in S}
\end{equation}
\end{theorem}
\begin{proof}
  Without loss of generality we can assume that $1 \in T$ (the coset
  representative of $H$ itself). By definition, $Hts = H\overline{ts}$
  which implies that $t s (\overline{ts})^{-1} \in H$ for all $t \in
  T, s \in S$. Hence, $S_H \subseteq H$ and $\gen{S_H} \leq H$, so the
  content of the statement lies in the other inclusion.
  
  Let $h \in H \leq G$ and observe that since $\gen{S} = G$ we have $h
  = s_1 s_2 \dotsm s_k$ where $s_i \in S \cup S^{-1}$ for $i = 1,
  \dotsc, k$. Define a sequence $t_1, t_2, \dotsc, t_{k + 1}$ of $k +
  1$ elements of $T$ as follows: $t_1 = 1$ and inductively $t_{i + 1} = \overline{t_i s_i}$. Furthermore, let $a_i = t_i s_i t_{i + 1}^{-1}$ for $i = 1, \dotsc, n$ and observe that
\begin{equation}
h = (t_1 s_1 t_2^{-1})(t_2 s_2 t_3^{-1}) \dotsm (t_n s_n t_{n + 1}^{-1}) t_{n + 1} = a_1 a_2 \dotsm a_n t_{n + 1}
\end{equation}
We now show that $a_i \in \gen{S_H}$ for $i = 1, \dotsc, n$, and that $t_{n + 1} = 1$, which
implies that $H \leq \gen{S_H}$. 

For each $i = 1, \dotsc, n$, either $s_i \in S$ or $s_i^{-1} \in S$.
In the first case we immediately get $a_i = t_i s_i (\overline{t_i
  s_i})^{-1} \in S_H$, and in the second case we have $H t_{i + 1}
s_i^{-1} = H \overline{t_i s_i} s_i^{-1} = H t_i s_i s_i^{-1} = H t_i$
which implies that $t_i = \overline{t_{i + 1} s_i^{-1}}$. Hence,
$a_i^{-1} = t_{i + 1} s_i^{-1} t_i^{-1} = t_{i + 1} s_i^{-1}
\left(\overline{t_{i + 1} s_i^{-1}}\right)^{-1} \in S_H$ and thus $a_i \in \gen{S_H}$.

Finally, since $h \in H$ and $\gen{S_H} \leq H$ we have $t_{n + 1} = (a_1 a_2 \dotsm a_n)^{-1} h \in H$, so $t_{n + 1}$ is the coset representative of $H$, and therefore $t_{n + 1} = 1$. Thus $H \leq \gen{S_H}$.
\end{proof}

Our situtation is that we have a group $G = \gen{S}$ acting on the finite set
$X$, and we want to use Schreier's Lemma to find the generators
(usually called the \emph{Schreier generators}) for the stabiliser
$G_{\alpha}$, where $\alpha \in X$. If we compute a Schreier tree for
the orbit $\alpha^G$ using \ref{alg:orbit} then we know that
\ref{alg:trace} can be used to find the transversal of the cosets of
$G_{\alpha}$ in $G$. 

For $p \in X$, let $t(p) \in G$ denote the result of \ref{alg:trace}.
Using the notation in Theorem \ref{schreiers_lemma} we then have $\bar{g} =
t(\alpha^g)$ and the transversal is $\set{t(p) \mid p \in \alpha^G}$,
so for $s \in S$ and $p \in \alpha^G$ the Schreier generator can be
expressed as
\begin{equation} 
t(p) s t(\alpha^{t(p) s})^{-1} = t(p) s t((\alpha^{t(p)})^s)^{-1} = t(p) s t(p^s)^{-1}
\end{equation}
and a generating set for $G_{\alpha}$ is
\begin{equation} \label{gen_stabiliser}
\set{t(p) s t(p^s)^{-1} \mid p \in \alpha^G, s \in S}
\end{equation}

\subsection{Computing a base and SGS}
Going back to our problem, if we have a partial base $B = (\alpha_1,
\dotsc, \alpha_n)$ for $G$ with corresponding partial strong
generating set $S$, then we can use Schreier's Lemma to solve our
problem. Using the notation from Theorem \ref{thm_leon}, we calculate
the Schreier generators for each $G^i$, using \eqref{gen_stabiliser} and add them to $S$, possibly
adding points to $B$, if some Schreier generator fixes the whole base,
in order to ensure that $B$ and $S$ are still partial. When this is
finished, we have $H^i = G^i$ for $i = 1, \dotsc, n$ and thus $B$ and
$S$ are complete. We can use \ref{alg:get_schreier_gen} to calculate a Schreier generator.

\begin{algorithm} 
\dontprintsemicolon
\caption{\texttt{GetSchreierGenerator}}
\SetKwFunction{OrbitElement}{OrbitElement}
\KwData{A group $G = \gen{S}$ acting on a finite set $X$, a Schreier tree $T$ for the orbit $\alpha^G$ of the point $\alpha \in X$, a $p \in X$ and a generator $s \in S$.}
\KwResult{The Schreier generator corresponding to $p$ and $s$}
\Begin
{
  $t_1 := \OrbitElement(T, p)$ \;
  $t_2 := \OrbitElement(T, p^s)$ \;
  \Return{$t_1 s t_2$}
}
\refstepcounter{algorithm}
\label{alg:get_schreier_gen}
\end{algorithm}

However, there is a problem with this simple approach, in that the
generating sets defined in \eqref{gen_stabiliser} can be very large,
and contain many redundant generators. A fraction of the Schreier
generators is usually enough to generate the stabiliser. In
\cite{hall59} it is shown that $[G : G_{\alpha}] - 1$ of the Schreier
generators for $G_{\alpha}$ are equal to the identity. For instance,
if for some point $p \in \alpha^G$ we have that $(p, p^s)$ is an edge
of the Schreier tree for $\alpha^G$, then the Schreier generator $t(p)
s t(p^s)^{-1}$ is the identity. Thus, the number of non-trivial
Schreier generators can be as large as $(\abs{S} - 1) [G : G_{\alpha}]
+ 1$, though some of them may be equal to each other.

In our case we calculate the Schreier generators for each $G^i$, using $S^{i - 1}$ in place of $S$, and since $S^{i - 1}$ are precisely the calculated Schreier generators for $G^{i - 1}$, the number of non-trivial Schreier generators for $G^i$ may be as large as 
\begin{equation} \label{num_schreier_gens}
1 + (\abs{S} - 1) \prod_{j = 0}^i \abs{\alpha_{j + 1}^{G^j}}
\end{equation}
Since the orbit sizes are only bounded by $\abs{X}$, we see that
\eqref{num_schreier_gens} is exponential in $\abs{X}$, and therefore
this method may not be efficient.

\subsection{Reducing the number of generators}
It is possible to reduce the number of generators at each step, so
that our generating set never grows too large. This can be done using
\ref{alg:reduce_schreier_gens}, which is due to Sims, and which can also be
found in \cite{sims98} and \cite{soicher98}.

\begin{algorithm} 
\dontprintsemicolon
\caption{\texttt{BoilSchreierGenerators}}
\KwData{A group $G$ acting on a finite set $X$, a partial base $B = (\alpha_1, \dotsc, \alpha_k)$ and corresponding partial strong generating set $S$ for $G$, an integer $1 \leq m \leq k$.}
\KwResult{A smaller partial strong generating set for $G$}
\Begin
{
  \For{$i := 1$ \KwTo $m$}
  {
    $T := S^{i - 1}$ \;
    \ForEach{$g \in T$}
    {
      \ForEach{$h \in T$}
      {
        \If{$\alpha_i^g = \alpha_i^h \neq \alpha_i$}
        {
          $S := (S \setminus \set{h}) \cup \set{g h^{-1}}$ \;
        }
      }
    }    
  }
  \Return{$S$}
}
\refstepcounter{algorithm}
\label{alg:reduce_schreier_gens}
\end{algorithm}

This algorithm reduces the generating set to size $\binom{\abs{X}}{2}
\in \OR{\abs{X}^2}$. With this algorithm, we can solve our problem using
\ref{alg:call_ss} and \ref{alg:pre_ss}. This, however, is not the
Schreier-Sims algorithm, which is a more clever and efficient method
of performing the same things.

\begin{algorithm} 
\dontprintsemicolon
\caption{\texttt{ComputeBSGS}}
\SetKwFunction{GetPartialBSGS}{GetPartialBSGS}
\SetKwFunction{Schreier}{Schreier}
\SetKwFunction{BoilSchreierGenerators}{BoilSchreierGenerators}
\SetKwData{SGS}{sgs}
\SetKwData{Base}{base}
\KwData{A group $G = \gen{S}$ acting on a finite set $X$.}
\KwResult{A base and strong generating set for $G$}.
\Begin
{
  $(\Base, \SGS) := \GetPartialBSGS(S, \emptyset)$ \;
  \For{$i := 1$ \KwTo $\abs{\Base}$}
  {
    $(\Base, \SGS) := \Schreier(\Base, \SGS, i)$ \;
    $\SGS := \BoilSchreierGenerators(\Base, \SGS, i)$ \;
  }
  \Return{$(\Base, \SGS)$}
}
\refstepcounter{algorithm}
\label{alg:call_ss}
\end{algorithm}

\begin{algorithm} 
\dontprintsemicolon
\caption{\texttt{Schreier}}
\SetKwFunction{GetSchreierGenerator}{GetSchreierGenerator}
\SetKwFunction{ComputeSchreierTree}{ComputeSchreierTree}
\SetKwFunction{NewBasePoint}{NewBasePoint}
\SetKwData{Tree}{tree}
\SetKwData{Gen}{gen}
\SetKwData{Point}{point}
\KwData{A group $G$ acting on a finite set $X$, a partial base $B = (\alpha_1, \dotsc, \alpha_k)$ and corresponding partial strong generating set $S$ for $G$, an integer $1 \leq i \leq k$ such that $G^j = H^j$ for $j = 0, \dotsc, i - 1$.}
\KwResult{Possibly extended partial base $B = (\alpha_1, \dotsc, \alpha_m)$ and corresponding partial strong generating $S$ set for $G$ such that $G^j = H^j$ for $j = 0, \dotsc, i$.}
\tcc{Assumes the existence of a function \texttt{NewBasePoint}$(g)$ that returns a point $p \in X$ such that $p^g \neq p$}
\Begin
{
  $T := S^{i - 1}$ \;
  $\Tree := \ComputeSchreierTree(T, \alpha_i)$ \;
  \ForEach{$p \in \alpha_i^{H^{i - 1}}$}
  {
    \ForEach{$s \in T$}
    {
      $\Gen := \GetSchreierGenerator(\Tree, p, s)$ \;
      \If{$\Gen \neq 1$}
      {
        $S := S \cup \set{\Gen, \Gen^{-1}}$ \;
        \If{$B^{\Gen} = B$}
        {
          $\Point := \NewBasePoint(\Gen)$ \;
          $B := B \cup \set{\Point}$ \;
        }
      }
    }    
  }
  \Return{$(B, S)$}
}
\refstepcounter{algorithm}
\label{alg:pre_ss}
\end{algorithm}

\section{Membership testing}
We now assume that we know a (complete) base $B = (\alpha_1, \dotsc,
\alpha_n)$ and a (complete) strong generating set $S$ for our group
$G$, and we present an efficient algorithm for determining if, given
an arbitrary group element $g$, it is true that $g \in G$. The
implicit assumption is of course that $G \leq H$ for some large group
$H$ and that $g \in H$, but since we are interested in matrix groups
this is always true with $H$ being some general linear group.

This algorithm is used in the Schreier-Sims algorithm together with
Theorem \ref{thm_leon}, as we shall see later. 

Recall that if we have a base then there is an associated stabiliser
chain, and as described in section \ref{section:bsgs}, if $g \in G$ we
can factorise $g$ as a product of coset representatives $g = u_n u_{n - 1}
\dotsm u_1$ where $u_i$ is the representative of $G^i g$ in $G^{i -
  1}$. Moreover, if we for each $i = 1, \dotsc, n$ have computed a
Schreier tree $T_i$ for $\alpha_i^{G^{i - 1}}$ then we can use
\ref{alg:trace} to compute each coset representative.

More specifically, if $g \in G$ then $g = g_1 t(\alpha_1^g)$ where, as
before, $t(p)$ is the output of \ref{alg:trace} on the point $p$ and
$g_1 \in G^1$. On the other hand, if $g \notin G$ then either
$\alpha_1^g \notin \alpha^G$ or $g_1 \notin G_1$. To test if $g \in G$
we can therefore proceed inductively, and first check if $\alpha_1^g
\in \alpha^G$ and if that is true then test whether $g_1 \in G_1$.
This is formalised in \ref{alg:membership}.

\begin{algorithm} 
\dontprintsemicolon
\caption{\texttt{Membership}}
\SetKwFunction{OrbitElement}{OrbitElement}
\SetKwData{Tree}{tree}
\SetKwData{Element}{element}
\SetKwData{Point}{point}
\KwData{A group $G$ acting on a finite set $X$, a base $B = (\alpha_1, \dotsc, \alpha_n)$, a Schreier tree $T_i$ for the orbit $\alpha_{i + 1}^{G^i}$ for each $i = 0, \dotsc, n - 1$, and a group element $g$.}
\KwResult{A residue $r$ and drop-out level $1 \leq l \leq n + 1$.}
\Begin
{
  $r := g$ \;
  \For{$i := 1$ \KwTo $n$}
  {
    \If{$\alpha_i^r \notin T_{i - 1}$}
    {
      \Return{$(r, i)$}
    }
    $\Element := \OrbitElement(T_{i - 1}, \alpha_i^r)$ \;
    $r := r \cdot \Element^{-1}$ \;
  }
  \Return{$(r, n + 1)$} \label{alg:membership_return}
}
\refstepcounter{algorithm}
\label{alg:membership}
\end{algorithm}

The terminology \emph{residue} and \emph{level} is introduced here,
with obvious meanings. As can be seen, the algorithm returns the level
at which it fails, as this is needed in the Schreier-Sims algorithm. Note that even if all $n$ levels are passed, it might happen that the residue $r \neq 1$ at line \ref{alg:membership_return}, which also indicates that $g \notin G$. 

In the literature, \ref{alg:membership} is usually referred to as
\emph{sifting} or \emph{stripping} of the group element $g$.

\section{The main algorithm}
Finally, we can now present the Schreier-Sims algorithm itself. It
uses a more efficient method of reducing the number of Schreier
generators considered, by making use of \ref{alg:membership}. 

Using
the notation from Theorem \ref{thm_leon}, we want to show that
$H^i_{\alpha_{i + 1}} = H^{i + 1}$ for each $i$, or equivalently that
all Schreier generators for $H^i$ are in $H^{i + 1}$. If we proceed
from $i = n - 1, \dotsc, 0$ instead of the other way, then for $H^n = 1$ we obviously already have a base and strong generating set,
so we can use \ref{alg:membership} to check if the Schreier generators
for $H^{n - 1}_{\alpha_n}$ are in $H^n$.

When we have checked all Schreier generators for $H^{n - 1}_{\alpha_n}$ we then have a
base and strong generating set for $H^{n - 1}$ by Theorem \ref{thm_leon},
since $H^{n - 1}_{\alpha_n} = H^n$, and we can therefore proceed
inductively downwards. This is shown in \ref{alg:ss_main} and \ref{alg:ss}.

\begin{algorithm} 
\dontprintsemicolon
\caption{\texttt{ComputeBSGS}}
\SetKwFunction{GetPartialBSGS}{GetPartialBSGS}
\SetKwFunction{SchreierSims}{SchreierSims}
\SetKwData{SGS}{sgs}
\SetKwData{Base}{base}
\KwData{A group $G = \gen{S}$ acting on a finite set $X$.}
\KwResult{A base and strong generating set for $G$}.
\Begin
{
  $(\Base, \SGS) := \GetPartialBSGS(S, \emptyset)$ \;
  \For{$i := \abs{\Base}$ \KwTo $1$}
  {
    $(\Base, \SGS) := \SchreierSims(\Base, \SGS, i)$ \;
  }
  \Return{$(\Base, \SGS)$}
}
\refstepcounter{algorithm}
\label{alg:ss_main}
\end{algorithm}

\begin{algorithm} 
\dontprintsemicolon
\caption{\texttt{SchreierSims}}
\SetKwFunction{GetSchreierGenerator}{GetSchreierGenerator}
\SetKwFunction{ComputeSchreierTree}{ComputeSchreierTree}
\SetKwFunction{NewBasePoint}{NewBasePoint}
\SetKwFunction{Membership}{Membership}
\SetKwFunction{SchreierSims}{SchreierSims}
\SetKwData{Tree}{tree}
\SetKwData{Gens}{gens}
\SetKwData{Gen}{gen}
\SetKwData{Point}{point}
\SetKwData{Residue}{residue}
\SetKwData{DropOut}{dropout}

\KwData{A group $G$ acting on a finite set $X$, a partial base $B =
  (\alpha_1, \dotsc, \alpha_k)$ and corresponding partial strong
  generating set $S$ for $G$, an integer $1 \leq i \leq k$ such that
  $H^{j - 1}_{\alpha_j} = H^j$ for $j = i + 1, \dotsc, k$, and
  Schreier trees $T^{j - 1}$ for $\alpha_j^{H^{j - 1}}$ for $j = i + 1, \dotsc,
  k$.}

\KwResult{Possibly extended partial base $B = (\alpha_1, \dotsc,
  \alpha_m)$ and corresponding partial strong generating $S$ set for
  $G$ such that $H^{j - 1}_{\alpha_j} = H^j$ for $j = i, \dotsc,
  m$ and Schreier trees $T^{j - 1}$ for $\alpha_j^{H^{j - 1}}$ for $j = i,
  \dotsc, m$.}

\tcc{Assumes the existence of a function \texttt{NewBasePoint}$(g)$ that returns a point $p \in X$ such that $p^g \neq p$}
\Begin
{
  $\Gens := S^i$ \;
  $T^{i - 1} := \ComputeSchreierTree(\Gens, \alpha_i)$ \; \label{alg:ss_comp_tree}
  \ForEach{$p \in \alpha_i^{H^{i - 1}}$}
  {
    \ForEach{$s \in \Gens$}
    {
      $\Gen := \GetSchreierGenerator(T^{i - 1}, p, s)$ \;
      \If{$\Gen \neq 1$}
      {
        $(\Residue, \DropOut) := \Membership(T^i, \dotsc, T^k, \Gen)$ \;
        \If{$\Residue \neq 1$}
        {
          $S := S \cup \set{\Gen, \Gen^{-1}}$ \;
          \If{$\DropOut = k + 1$}
          {
            $\Point := \NewBasePoint(\Gen)$ \;
            $B := B \cup \set{\Point}$ \;
          }
          \For{$j := r$ \KwTo $i + 1$}
          {
            $(B, S) := \SchreierSims(B, S, j)$ \;
          }
        }
      }    
    }
  }
  \Return{$(B, S)$}
}
\refstepcounter{algorithm}
\label{alg:ss}
\end{algorithm}

\subsection{Matrix groups}
As can be seen in the given algorithms, the main part of the
Schreier-Sims algorithm is indepedent of the particular type of group,
but we are interested in the situation where $G \leq \GL(d, q)$ for
some $d \geq 1$ and some $q = p^r$ where $p$ is a prime number and $r
\geq 1$. Here, $d$ is called the \emph{degree} of $G$ and is the
number of rows (columns) of the matrices in $G$, and $q$ is the finite
field size, ie $G$ contains matrices over $\GF(q) = \F_q$. In this
setting, $G$ acts faithfully from the right on the vector space $X = \F_q^d$, by multiplication of a row vector with a matrix. We denote the standard base in $\F_q^d$ by $\mathbf{e}_1, \mathbf{e}_2, \dotsc, \mathbf{e}_d$.

What remains to specify is the algorithm \texttt{NewBasePoint}, which
depends on the particular point set $X$ and the action that is used.
To construct this algorithm, we observe that given a matrix $M \neq I
\in G$, if $M_{ij} \neq 0$ for some $i \neq j$ then $\mathbf{e}_i M \neq \mathbf{e}_i$
so the row vector $\mathbf{e}_i$ is a point moved by $M$. If $M$ instead is
diagonal, but not a scalar matrix, then $M_{ii} \neq M_{jj}$ for some
$i \neq j$, and $(\mathbf{e}_i + \mathbf{e}_j) M \neq \mathbf{e}_i + \mathbf{e}_j$ so this row vector is
moved by $M$. Finally if $M$ is scalar, then $M_{11} \neq 1$ since $M
\neq I$ and thus $\mathbf{e}_1 M \neq \mathbf{e}_1$. This translates directly into \ref{alg:new_base_point}.

\begin{algorithm} 
\dontprintsemicolon
\caption{\texttt{NewBasePoint}}
\KwData{A matrix $M \neq I \in G \leq \GL(d, q)$.}
\KwResult{A row vector $\mathbf{v} \in \F_q^d$ such that $\mathbf{v} M \neq \mathbf{v}$}.
\SetKw{KwAnd}{and}
\Begin
{
  \For{$i := 1$ \KwTo $d$}
  {
    \For{$j := 1$ \KwTo $d$}
    {
      \If{$i \neq j$ \KwAnd $M_{ij} \neq 0$}
      {
        \Return{$\mathbf{e}_i$}
      }
    }
  }
  \For{$i := 1$ \KwTo $d$}
  {
    \For{$j := 1$ \KwTo $d$}
    {
      \If{$i \neq j$ \KwAnd $M_{ii} \neq M_{jj}$}
      {
        \Return{$\mathbf{e}_i + \mathbf{e}_j$}
      }
    }
  }
  \Return{$\mathbf{e}_1$}
}
\refstepcounter{algorithm}
\label{alg:new_base_point}
\end{algorithm}

\section{Complexity analysis}
We now analyse the time complexity of some of the given algorithms, but we
will not be interested in space complexity. For the analysis we assume
that the external functions \texttt{Tree}, \texttt{AddChild} and
\texttt{EdgeLabel}, that depend on the particular datastructure used,
take $\Theta(1)$ time. This is a reasonable assumption and it is satisfied
in the code for the project, which uses hash tables to implement the Schreier trees. We also assume that the datastructure
used for sets is a sorted list, so that elements can be found, added and
removed in logarithmic time using binary search. This assumption is
satisfied in \GAP.

Moreover, since we are working in a matrix group $G = \gen{S} \leq
\GL(d, q)$ that acts on the vector space $X = \F_q^d$, we know that
the multiplication of two group elements takes $\Theta(d^3)$ time and
that the action of a group element on a point takes 
$\Theta(d^2)$ time, under the assumption that the multiplication of two
elements from $\F_q$ takes $\Theta(1)$ time. This assumption is
reasonable since we are mostly interested in the case when $q$ is not too large, so that one field element can be stored in a machine word and manipulated in constant time. Then
we also see that testing equality between two group elements takes
$\OR{d^2}$ time and testing equality between two points takes 
$\OR{d}$ time.

Consider first \ref{alg:orbit}. We know that it is a breadth-first
search and we know that a simple breadth-first search on a graph
$\mathcal{G} = (V, E)$ takes $\OR{\abs{V} + \abs{E}}$ time. In our case
we have $\abs{E} = \abs{V} \abs{S} = \abs{X} \abs{S}$, so the edges
are dominating. We also see that line \ref{alg:orbit_innerloop1}
takes $\Theta(d^2)$ time and line \ref{alg:orbit_innerloop2} takes
$\OR{\log{\abs{S}}}$ time since the size of the set
\textsf{children} is bounded by $\abs{S}$. The former line is therefore dominating, and thus we have that
\ref{alg:orbit} takes $\OR{\abs{X}\abs{S} d^2}$ time.

For \ref{alg:trace} we just note that in a tree with vertex set $V$,
the depth of any node is bounded by $\abs{V}$. In our case the
vertices of the Schreier tree is the points of the orbit, which may be
the whole of $X$ if the action is transitive. Therefore we see that
\ref{alg:trace} takes time $\OR{\abs{X} d^4}$.

It is obvious that \ref{alg:new_base_point} takes time $\OR{d^3}$, and
it is also obvious that \ref{alg:get_schreier_gen} takes time
$\OR{\abs{X} d^2 + d^3} = \OR{\abs{X} d^2}$. In \ref{alg:partial_bsgs}
we see that line \ref{alg:partial_bsgs_innerloop2} takes time $\OR{1}$
because the base is just a list and can therefore be augmented in
constant time. We see that line \ref{alg:partial_bsgs_innerloop1}
takes time $\OR{\abs{B} d^3}$ since we might have to multiply every
base point with the group element and check equality, in the case
where it actually fixes the base. Thus, \ref{alg:partial_bsgs} takes
time $\OR{\abs{S} \log{\abs{S}} \abs{B} d^3}$.

Since \ref{alg:pre_ss} is not used in the \GAP~package of the project,
we skip the complexity analysis of it and its related algorithms. We also skip any detailed complexity
of the Schreier-Sims algorithm itself, since we will see that for
matrix groups, this is futile anyway. The case of permutation groups is
analysed in \cite{butler91} and the result is that the Schreier-Sims
algorithm has time complexity $\OR{\abs{X}^5}$. 

The essential problem
when using the algorithm for matrix groups is then that $X = \F_q^d$
so $\abs{X} = q^d$ and the Schreier-Sims algorithm is thus exponential
in $d$. This is quite disheartening, since it is the case when $d$
grows large that we are interested in, more often than when $q$ is
large. However, it is possible to make the implementation fast enough for many practical purposes, so all is not lost.

\chapter{The probabilistic Schreier-Sims algorithm}
It is well-known in computer science that probabilistic algorithms
often turn out to be much simpler in structure than the corresponding
deterministic algorithms for the same problem, while still producing
good solutions. That the algoritms are simpler usually means that they
can be made more efficient and are easier to implement. 

On the other hand, the drawback with the probabilistic approach is
that there is a non-zero probability that the algorithms in some sense
can ''lie'', ie return incorrect solutions. An introduction to the
subject and the various complexity classes can be found in
\cite{papadimitriou94}.

The probabilistic Schreier-Sims algorithm is a good example of this
situation. It was first described in \cite{leon80}, and the idea comes from the following
\begin{theorem} \label{thm_rand_schreier1}
  Let $G$ be a group acting on the finite set $X$, and let $B =
  (\alpha_1, \dotsc, \alpha_k)$ be a partial base and $S$ a partial
  strong generating set for $G$. For $i = 1, \dotsc, n$, let $G^i =
  G_{\alpha_1, \dotsc, \alpha_i}$, $S^i = S \cap G^i$, $H^i =
  \gen{S^i}$, let $G = G^0 = H^0$ and let $T^i$ be a right transversal for the cosets of
  $H^{i - 1}_{\alpha_i}$ in $H^{i - 1}$. Then $\abs{\prod_{i = 1}^n T^i}$ divides $\abs{G}$.
\end{theorem}
\begin{proof}
Since $\abs{T^i} = [H^{i - 1} : H^{i - 1}_{\alpha_i}]$ we have 
\begin{equation}
\begin{split}
\abs{G} &= \prod_{i = 1}^n [H^{i - 1} : H^i] = \prod_{i = 1}^n [H^{i - 1} : H^{i - 1}_{\alpha_i}] [H^{i - 1}_{\alpha_i} : H^i] = \\
&= \left( \prod_{i = 1}^n \abs{T^i} \right) \left( \prod_{i = 1}^n [H^{i - 1}_{\alpha_i} : H^i] \right)
\end{split}
\end{equation}
and thus the theorem follows.
\end{proof}

\begin{cl} \label{thm_random_schreier2}
  Let $G$ be a group acting on the finite set $X$, and let $B =
  (\alpha_1, \dotsc, \alpha_k)$ be a partial base and $S$ a partial
  strong generating set for $G$. If $B$ and $S$ are not complete and $g \in G$ is a uniformly random element, then the probability that \ref{alg:membership} returns residue $r \neq 1$ when given $g$ is at least $1/2$.
\end{cl}
\begin{proof}
  We see that even if $B$ and $S$ are not complete, we can compute
  Schreier trees for the orbits $\alpha_{i + 1}^{H^i}$ for $i = 0,
  \dotsc, n - 1$ and feed them to \ref{alg:membership} which then
  tries to express $g$ in terms of the corresponding right
  transversals $T^i$ for $[H^{i - 1} : H^{i - 1}_{\alpha_i}]$.
  
  Thus, \ref{alg:membership} checks if $g \in \prod_{i = 1}^n T^i$ and
  by Theorem \ref{thm_rand_schreier1}, if $\prod_{i = 1}^n T^i$ is not
  the whole $G$, then it contains at most half of the elements of $G$.
  Therefore, since $g$ is uniformly random, we have $\Pr{g \notin
    \prod_{i = 1}^n T^i} \geq 1/2$.
\end{proof}

This suggests a probabilistic algorithm for computing a base and a
strong generating set. Given a partial base and partial strong
generating set for $G$, compute Schreier trees and chose random
elements uniformly from $G$. Use \ref{alg:membership} on each element,
and if it returns a non-trivial residue, add it to the partial SGS,
and possibly augment the base.

If the base and SGS are complete, then of course the residue will be trivial. On the other hand, if the base and SGS are not complete, then by
Corollary \ref{thm_random_schreier2}, the probability that $k$
consecutive random elements have trivial residues is less than
$2^{-k}$. We can thus choose a large enough $k$ for our purposes and
assume that the base and SGS are complete when we have stripped $k$
consecutive random elements to the identity. This is formalised in
\ref{alg:random_schreier}.

\begin{algorithm} 
\dontprintsemicolon
\caption{\texttt{RandomSchreierSims}}
\SetKwFunction{ComputeSchreierTree}{ComputeSchreierTree}
\SetKwFunction{NewBasePoint}{NewBasePoint}
\SetKwFunction{Random}{Random}
\SetKwFunction{Membership}{Membership}
\SetKwData{Tree}{tree}
\SetKwData{Sifts}{sifts}
\SetKwData{Element}{element}
\SetKwData{Point}{point}
\SetKwData{Residue}{residue}
\SetKwData{DropOut}{dropout}

\KwData{A group $G$ acting on a finite set $X$, a partial base $B =
  (\alpha_1, \dotsc, \alpha_k)$ and corresponding partial strong
  generating set $S$ for $G$, an integer $m \geq 1$ and Schreier trees $T^{j - 1}$ for $\alpha_i^{H^{i - 1}}$ for $i = 1, \dotsc, k$.}

\KwResult{Possibly extended partial base $B = (\alpha_1, \dotsc,
  \alpha_n)$ and corresponding partial strong generating $S$ set for
  $G$ such that $m$ consecutive random elements have been stripped to identity with respect to $B$ and $S$.}

\tcc{Assumes the existence of a function \texttt{NewBasePoint}$(g)$ that returns a point $p \in X$ such that $p^g \neq p$.}
\tcc{Assumes the existence of a function \texttt{Random}$(G)$ that returns a uniformly random element from the group $G$.}
\Begin
{
  $\Sifts := 0$ \;
  \While{$\Sifts < m$}
  {
    $\Element := \Random(G)$ \;
    $(\Residue, \DropOut) := \Membership(T^1, \dotsc, T^k, \Element)$ \;
    \eIf{$\Residue \neq 1$}
    {
      $S := S \cup \set{\Element, \Element^{-1}}$ \;
      \If{$\DropOut = k + 1$}
      {
        $\Point := \NewBasePoint(\Element)$ \;
        $B := B \cup \set{\Point}$ \;
        $k := k + 1$ \;
      }
      
      \For{$i := 1$ \KwTo $k$}
      {
        $T^i := \ComputeSchreierTree(S^i, \alpha_i)$ \;
      }
      $\Sifts := 0$ \;
    }
    {
      $\Sifts := \Sifts + 1$ \;
    }
  }
  \Return{$(B, S)$}
}
\refstepcounter{algorithm}
\label{alg:random_schreier}
\end{algorithm}

Indeed, \ref{alg:random_schreier} is simpler than \ref{alg:ss} and
implementations are usually faster. However, it is not clear how the
algorithm \texttt{Random} should be constructed. First of all, we need
random bits, and the generation of pseudo-random bits is an old
problem in computer science, see \cite{papadimitriou94} and \cite{knuth97}, which we will not dwell into. Instead, we will assume that uniformly random bits are available.

\section{Random group elements}

The generation of random group elements is also difficult in general,
and an introduction to the topic can be found in \cite{babai96}. If
one knows a base and strong generating set one can easily generate
random elements by using the factorisation described in section
\ref{section:bsgs}, select random elements from each transversal and
multiply them. In our case we do not yet have a base and an SGS, since that is what we are trying to compute, and
then the only known algorithm that can \emph{provably} generate (almost)
uniformly random group elements is described in \cite{babai91} and
runs in time $\OR{d^{10} (\log{q})^5}$ for subgroups of $\GL(d, q)$, and
this is too slow for our purposes.

The practical algorithms are instead \emph{heuristics}, see
\cite{kann99}, for which there are no \emph{proven} guarantees for any
uniformly randomness, but instead experimentation and statistical
analysis has shown them to be good. It appears that the most
successful algorithm is the \emph{product replacement} algorithm, also
called ''Shake'', which is originally an idea by Charles Leedham-Green and
Leonard Soicher, and is described in \cite{lg95}. 

We will not dwell deeper into the analysis of this algorithm, but
the algorithm itself is quite simple. Given a group $G =
\gen{g_1, \dotsc, g_n}$, the ''Shake'' algorithm maintains a global
variable $S = (a_1, \dotsc, a_m) \in G^m$ for some $m \geq n$, and
then each call to \ref{alg:shake} returns a proposed random group
element.
\begin{algorithm} 
\dontprintsemicolon
\caption{\texttt{Shake}}
\SetKwFunction{RandomInteger}{RandomInteger}
\SetKwData{State}{state}

\KwData{The global state $S = (a_1, \dotsc, a_m) \in G^m$}

\KwResult{An element in $G$ that hopefully is a good approximation of a uniformly random element from $G$.}

\tcc{Assumes the existence of a function \texttt{RandomInteger}($k$) that returns a uniformly random integer from the set $\set{0, \dotsc, k}$}
\Begin
{
  $i := \RandomInteger(m)$ \;
  \Repeat{$i \neq j$}
  {
    $j := \RandomInteger(m)$ \;
  }
  \eIf{$\RandomInteger(1) = 0$}
  {
    $b := a_i a_j$ \;
  }
  {
    $b  := a_j a_i$ \;
  }
  $a_i := b$ \;
  \Return{$b$}
}
\refstepcounter{algorithm}
\label{alg:shake}
\end{algorithm}

Evidently, the algorithm is very cheap in terms of time and space. The
questions that remain are how the state $S$ should be initialised and
how $m$ should be chosen, and the ability of the algorithm to generate
uniformly random group elements are highly dependent on the answers.
In \cite{lg95} the suggestion is that $m = \max (10, 2n + 1)$ and that
$S$ is initialised to contain the generators $g_1, \dotsc, g_n$, the
rest being identity elements. Moreover, the state should be
initiliased by calling \ref{alg:shake} a number $K$ of times and
discarding the results, and the suggestion is that $K \geq 60$.

In \cite{niemeyer01} this algorithm is also described, as well as a
variant of it, called ''Rattle'', which is due to Leedham-Green. The
algorithms are compared and ''Rattle'' is found to be slightly better,
but its running time is a bit higher. In \GAP, there is an
implementation of the ''Shake'' algorithm, which is used in this
project. An implementation of ''Rattle'' was also made, but it turned
out that it was not efficient enough.

\chapter{The Schreier-Todd-Coxeter-Sims algorithm}

Coset enumeration is one of the oldest algorithms in group theory, and
it was first described in \cite{tc36}. Given a finitely presented
group $G$ and $H \leq G$, the aim of coset enumeration is to construct
the list of the cosets of $H$ in $G$, ie to find a set $X$ with a
point $p \in X$ such that $G$ acts transitively on $X$ and such that
$H = G_p$. The set $X$ is usually taken to be $\set{1, \dotsc, n}$ for
some $n$ and $p = 1$. If the index of $H$ is not finite, then the algorithm will usually not terminate.

The coset enumeration is a trial-and-error process, and quite a few
strategies have been developed over the years. We will not be
concerned about the actual algorithm, since it is beyond the scope of this report, and there is a good
implementation in \GAP, which is used in this project. A more detailed
description of coset enumeration can be found in \cite{sims98}.

However, it is important to be aware of the overall structure of coset
enumeration. It works by letting the elements of $G$ act on $X$ and
introducing new cosets when necessary. Then it may happen that some
cosets turn out to be the same and thus are identified in the list.
The structure of the algorithm makes it possible to exit prematurely
when a given number of cosets have been defined, and this can be used
if one knows an upper bound on the number of cosets.

In our case, we want to use the last case of Theorem \ref{thm_leon}.
At level $i$ of Schreier-Sims algorithm, we know that the base and SGS
are complete for higher levels, and want to verify that $H^{i + 1} =
H^i_{\alpha_{i + 1}}$, which by the theorem is equivalent to $[H^i :
H^{i + 1}] = \abs{\alpha_{i + 1}^{H^i}}$ and since we can compute a
Schreier tree for $\alpha_{i + 1}^{H^i}$ we know the orbit size. Thus,
if we can compute $[H^i : H^{i + 1}]$ using coset enumeration and it turns out to be equal the
orbit size, it is unnecessary to compute any Schreier generators and
check for membership.

All this is formalised in \ref{alg:stcs} where we use the option to
exit the coset enumeration after $M \abs{\alpha_{i + 1}^{H^i}}$ cosets have
been defined, for some rational number $M \geq 1$, since this
certainly is an upper bound. In our implementation, the number $M$ can
be specified by the user, and the default is $M = 6/5$. This
value comes from \cite{leon80}, where the algorithm was introduced and
some experimentation was carried out to determine a good value of $M$.

\begin{algorithm} 
\dontprintsemicolon
\caption{\texttt{SchreierToddCoxeterSims}}
\SetKwFunction{GetSchreierGenerator}{GetSchreierGenerator}
\SetKwFunction{ComputeSchreierTree}{ComputeSchreierTree}
\SetKwFunction{NewBasePoint}{NewBasePoint}
\SetKwFunction{Membership}{Membership}
\SetKwFunction{ToddCoxeter}{ToddCoxeter}
\SetKwFunction{SchreierToddCoxeterSims}{SchreierToddCoxeterSims}
\SetKwData{Tree}{tree}
\SetKwData{Gens}{gens}
\SetKwData{Gen}{gen}
\SetKwData{Point}{point}
\SetKwData{Table}{table}
\SetKwData{Residue}{residue}
\SetKwData{DropOut}{dropout}

\KwData{A group $G$ acting on a finite set $X$, a partial base $B =
  (\alpha_1, \dotsc, \alpha_k)$ and corresponding partial strong
  generating set $S$ for $G$, an integer $1 \leq i \leq k$ such that
  $H^{j - 1}_{\alpha_j} = H^j$ for $j = i + 1, \dotsc, k$, and
  Schreier trees $T^{j - 1}$ for $\alpha_j^{H^{j - 1}}$ for $j = i + 1, \dotsc,
  k$.}

\KwResult{Possibly extended partial base $B = (\alpha_1, \dotsc,
  \alpha_m)$ and corresponding partial strong generating $S$ set for
  $G$ such that $H^{j - 1}_{\alpha_j} = H^j$ for $j = i, \dotsc,
  m$ and Schreier trees $T^{j - 1}$ for $\alpha_j^{H^{j - 1}}$ for $j = i,
  \dotsc, m$.}

\tcc{Assumes the existence of a function \texttt{NewBasePoint}$(g)$ that returns a point $p \in X$ such that $p^g \neq p$.}
\tcc{Assumes the existence of a function \texttt{ToddCoxeter}$(U_1, U_2, k)$ that performs coset enumeration on $G = \gen{U_1}$ and $H = \gen{U_2} \leq G$, exiting when $k$ cosets have been defined.}
\Begin
{
  $\Gens := S^i$ \;
  $T^{i - 1} := \ComputeSchreierTree(\Gens, \alpha_i)$ \;
  \ForEach{$p \in \alpha_i^{H^{i - 1}}$}
  {
    \ForEach{$s \in \Gens$}
    {
      $\Table := \ToddCoxeter(S^i, S^{i + 1}, \abs{T^i} + 1)$ \;
      \If{$\abs{\Table} = \abs{T^i}$}
      {
        \Return{$(B, S)$}
      }
      $\Gen := \GetSchreierGenerator(T^{i - 1}, p, s)$ \;
      \If{$\Gen \neq 1$}
      {
        $(\Residue, \DropOut) := \Membership(T^i, \dotsc, T^k, \Gen)$ \;
        \If{$\Residue \neq 1$}
        {
          $S := S \cup \set{\Gen, \Gen^{-1}}$ \;
          \If{$\DropOut = k + 1$}
          {
            $\Point := \NewBasePoint(\Gen)$ \;
            $B := B \cup \set{\Point}$ \;
          }
          \For{$j := r$ \KwTo $i + 1$}
          {
            $(B, S) := \SchreierToddCoxeterSims(B, S, j)$ \;
          }
        }
      }    
    }
  }
  \Return{$(B, S)$}
}
\refstepcounter{algorithm}
\label{alg:stcs}
\end{algorithm}

The Schreier-Todd-Coxeter-Sims algorithm is known to perform
particularly good when the input partial base and SGS are actually
already complete. Thus, it can be used for \emph{verification} of the output
from a probabilistic algorithm, and in this project it is mainly used
in that way, to verify the output from the probabilistic algorithm
described earlier.

\chapter{Implementation and optimisation}

We will now consider some more practical issues regarding the actual
implementation and optimisation of the code in the project.

In \GAP, a very fast implementation of the Schreier-Sims algorithm already
exists. It is described in \cite{seress03} and is a heuristic based on
the algorithm described in \cite{seress91}. As with many other
algorithms in \GAP, it works only with permutation groups, so for a
general group $G$ it first calculates a faithful permutation
representation, ie an injective homomorphism $\lambda : G \to
\Sym{X}$, for some finite set $X$, which exists due to the well-known theorem by Cayley. It
then work with the image $\lambda(G)$, or rather with $H \leq S_{\abs{X}}$ such that $H \cong \lambda(G)$.

Maybe the most important motivation for the \GAP~package developed in
this project is the idea that if $G$ is a matrix group, and one work
with the matrices directly instead of first converting to
permutations, then perhaps one could use the additional structure, that
is otherwise thrown away, and come up with a faster algorithm. The code is 
therefore created with the basic assumption that all group elements are
invertible matrices over a finite field, and it uses the internal
\GAP~representation of such matrices.

\section{Schreier tree creation}
An important issue is how the Schreier trees are managed. Actually,
what we really want is the transversals defined by the Schreier trees,
so a possibility is to store not a tree, but just a list of all points
and for each point the element that moves the root to this point. This
can be realised as a tree with height 1, where the edge labels are not
the generators, as in the case of a Schreier tree, but the group
elements, so in general this will take up more memory than the
Schreier tree. On the other hand, to find a coset representative in a
Schreier tree we need to follow a path from a point to the root, which
takes logarithmic time on average, and may take linear time if the
tree is not properly balanced, but in the case where we store the
coset representatives as edge labels, then of course it takes constant
time to find them later. 

To conclude, we have a trade-off between time and space, and in this
project, both of the above strategies are available to the user of the
package. The concrete representation of the trees are as hash tables,
where the keys are the points and the values are the edge labels of
the unique edges directed towards the root, ie we use
''back-pointers''. This makes it possible to perform in constant time
the common task of checking if a given point is in the orbit defined
by a Schreier tree, since this is just a hash table lookup.

\subsection{Extending vs creating}

If we consider \ref{alg:ss} in more detail, then we notice that the
sets $S^i$ are only augmented if they are changed. Therefore, at line
\ref{alg:ss_comp_tree} we could save the tree that is computed and, at
the next time we arrive there, only extend the tree using the new
generators, if there are any.

This may or may not give a faster algorithm. Evidently, there will be
less work in \ref{alg:orbit} to extend a Schreier tree than to create
a new one, but it may happen that the tree becomes more balanced if it
is recomputed using all generators than if it is extended. If the tree
is not balanced, then \ref{alg:trace} will take more time, and this is
in fact the function where most time is spent. In \cite{seress03} it
is claimed that this problem with unbalanced trees are more common for
matrix groups, and empirical studies in this project has shown that in
most cases, it is better to recompute the Schreier trees every time.

\subsection{Shallow trees}
There are algorithms for the creation of Schreier trees that are
guaranteed to make the Schreier trees \emph{shallow}, ie balanced, so
that the they have worst-case logarithmic height. This is crucial if
one wants to have good worst-case complexity, and two algorithms are
described in \cite{seress03}, one deterministic and one probabilistic.
In this project the deterministic algorithm have been implemented,
which is also described in \cite{seress91}. It is too complicated to
be included here in more detail, but the essential idea is to choose a
different set of edge labels, rather than the given generators.

\section{Orbit sizes}
As have been noted earlier, using the Schreier-Sims algorithm for a
matrix group $G \leq \GL(d, q)$ is in a sense doomed from the
beginning, since the complexity is exponential in $d$ when $G$ acts on
$\F_q^d$. One manifestation of doom in this case is that the orbits
may become huge, something that does not happen for permutation
groups. This will make our Schreier trees huge and a large amount of Schreier generators must be created, so the algorithm will be slow.

\subsection{Alternating actions}
To avoid large orbits, one can use another action of $G$. However,
only if the action is faithful is it a permutation representation of
$G$, and this is needed if the Schreier-Sims algoritm is going to
work, otherwise we will calculate a base and strong generating set for the quotient of $G$ with the kernel of the action.

In \cite{butler76}, a clever trick was introduced where base points are
chosen alternatingly as one-dimensional subspaces (lines) and vectors
from those lines. When a vector $v = (v_1, \dotsc, v_d) \in \F_q^d$ is
chosen as base point, it is preceded in the base by the line
$\gen{v}$. The action of $G$ on the lines is known as the
\emph{projective action} and of course it is not faithful, but since the next
point $v$ is from the kernel this does not matter. Also, the subspace
$\gen{v}$ has a canonical representative $(1, v_2 v_1^{-1}, \dotsc,
v_d v_1^{-1})$, which is trivial to find from $v$. Therefore, in terms
of time, the projective action is not particularly more expensive than the action on points.

Now, if $k = \abs{v^G}$ then $\abs{\gen{v}^G} \abs{v^{G_{\gen{v}}}} = k$ and $G_{\gen{v}}$
is the stabiliser of the line containing $v$. This implies that if $u \in
v^{G_{\gen{v}}}$ then $u$ is also on that line, so $u = mv$ where $m
\in \F_q$. We see that $M = \set{m \in \F_q \mid mv \in v^{G_{\gen{v}}}}$ is a
subgroup of $\F_q^{*}$ and thus $l = \abs{v^{G_{\gen{v}}}} = \abs{M}$ divides $q - 1$.

Instead of one orbit of size $k$ we therefore have two orbits of size
$k/l$ and $l$, and since $k \geq k/l + l$ whenever $l \geq 2$ we will for
example need to compute and sift fewer Schreier generators. On the
other hand, our base will probably be longer when using this trick,
and empirical studies in this project has shown that it is not always
a good idea, but it nevertheless implemented and can be used.

\subsection{Eigenspaces}
Another method for producing smaller orbits is described in \cite{murray95}. The idea is to choose the base points to be eigenvectors of the given generator matrices.

Recall from linear algebra that the \emph{characteristic polynomial}
of a matrix $A \in \GL (d, q)$ is $c_A(x) = \det (xI - A)$ where $I$
is the identity matrix, and an \emph{eigenvalue} of $A$ is a root of
$c_A(x)$.  Normally, and \emph{eigenvector} is defined as a vector $v
\in \F_q^d$ such that $v^A = \lambda v$ where $\lambda$ is an
eigenvalue of $A$, or equivalently $v^{g(A)} = 0$ where $g(x)$ is a
linear factor of $c_A(x)$. Here we use a more general version of the
latter definition, where we allow factors of higher degree as well as
linear factors of $c_A(x)$.

We see that for a linear factor $g(x) = x - \lambda$ and a
corresponding eigenvector $v$, the size of the orbit $v^{\gen{A}}$ is
a divisor of $q - 1$, since if $u \in v^{\gen{A}}$ then $u = v^{\alpha
  A} = \alpha \lambda v$, and as before the possible values of
$\alpha$ form a subgroup of $\F_q^{*}$. More generally, for a factor
of $c_A(x)$ of degree $m$, the size of the orbit $v^{\gen{A}}$ for an
eigenvector $v$ is bounded above by $q^m - 1$. To get the smallest
possible orbits, we should therefore choose eigenvectors corresponding
to factors of as small degree as possible, and linear factors are easy
to find since they correspond to eigenvalues, which are easy to
compute.

Note that we have only considered the orbits of groups generated by
single matrices, and the orbit of a matrix group $G = \gen{S}$ need
not be small just because we choose as base point an eigenvector of
one of the generators in $S$, but if we choose a base point that is an
eigenvector of several generators, then it turns out that the orbit
size are more often small. In \cite{murray95} these issues are
investigated and experimented with in some detail, and a heuristic for
finding base points that will hopefully give small orbits is
developed. This has also been implemented in \GAP, and in this project
it is possible to use that algorithm. However, it is not always a good
idea to use it, since there is some overhead, and it is not certain
that the orbit sizes will actually be smaller.

\section{Further developments}

The package contains implementations of two more advanced algorithms,
the so-called \emph{Verify} routine by Sims, which is an algorithm for
verifying if a proposed base and SGS are complete, and the nearly
linear time algorithm for finding a base and an SGS. The complete
descriptions of these algorithms are beyond the scope of this report,
but for completeness we include brief accounts on them.

\subsection{The Verify routine}

This algorithm is due to Sims, and has never been published, but it is
described in \cite{seress03}. Given a group $G = \gen{S}$ acting on
the finite set $X$, a point $\alpha \in X$ and a subgroup $H = \gen{S^{\prime}} \leq G_{\alpha}$, it
checks whether $H = G_{\alpha}$. If this is not the case, the
algorithm computes $g \in G_{\alpha} \setminus H$ to witness that. To
check a whole proposed base and SGS, one then uses the third case of
Theorem \ref{thm_leon} and checks each level.

The algorithm is quite involved, both theoretically and when it comes
to implementing it. It has a recursive nature, inducting on $S \setminus
S^{\prime}$, and it involves things like changing base points and
computing block systems.

\subsection{The nearly linear time algorithm}
For permutation groups $G$ acting on $X$, there is an algorithm for
computing a base and an SGS that runs in nearly linear time, ie linear
in $\abs{X}$ except for some logarithmic factor of $\abs{X}$ and
$\abs{G}$. It is the best known algorithm in terms of time complexity,
and is described in \cite{seress91} as well as in \cite{seress03} and it
is probabilistic. The complexity is achieved by using shallow
Schreier trees, a fast probabilistic verification algorithm, and by
sifting only a few randomly selected Schreier generators. Otherwise,
the algorithm is quite similar to the algorithms we have described
earlier for computing a base and an SGS.

The random Schreier generators are computed using \emph{random
  subproducts} which are described in \cite{seress03} as well as by
selecting random group elements, and the algorithm relies on some
theorems relating the number of Schreier generators to compute to the
given probability of correctness.

\chapter{Performance}

As observed earlier, one of the objectives of the project was to make
an implementation that hopefully would be faster than the one already
existing in \GAP. To determine if this objective was met, the
implementation has been benchmarked and compared with the built-in
\GAP~implementation.

The algorithm has been used to compute a base and SGS for some matrix
groups that are easy to construct in \GAP, and the generating sets
that were used were the standard generating sets from the
\GAP~library. The main test groups were classical groups: the general
and special linear groups $\GL(d, q)$ and $\SL(d, q)$ and the general
and special orthogonal groups $\GO(d, q)$ and $\SO(d, q)$, for various
(small) $d$ and $q$. The algorithm was also tested on some Suzuki
groups $\Sz(q)$ (where $q$ is a non-square power of $2$) and some Ree
groups $\Ree(q)$ (where $q = 3^{1 + 2m}$ for some $m > 0$).

Another, perhaps more realistic, test of the performance was done by
running the algorithm on randomly formed sets of invertible matrices
of a given dimension over a given finite field.

The benchmarks was carried out on two quite standard PC computers, the first test on a computer with an AMD Athlon
CPU running at $2$ GHz and with $1$ GB of physical RAM, and the second test on a computer with an Intel Pentium 4 CPU running at $2.8$ GHz and also with $1$ GB of physical RAM. It is likely
that these were the only important parameters, since \GAP~is mainly
CPU intensive, and swapping was avoided duringd the benchmark, so the hard disk speed should be a negligible factor. The
\GAP~installation tests reported GAP4stones values of $194624$ and $253581$, respectively.

The details of the benchmarks is shown in the appendix. During the first test, the existing algorithm was the faster one for most groups, but for some groups
our implementation was faster, most notably for the Suzuki groups. It
should be noted that all running times for our implementation comes
from using the same algorithm options, and it is possible to get
better times, especially for the smaller groups, by elaborating with
other options. The second test also indicated that the existing implementation is faster in most cases.

In terms of memory, no rigorous benchmark has been performed, but
during the above benchmark some simple checks of the amount of memory
allocated by the \GAP~process were performed. It seemed that our
implementation consumed less memory overall, and in some cases the
difference was large. Indeed, since we wanted to avoid swapping during
the benchmark, no larger Suzuki group than $\Sz(32)$ could be checked during the first test, since the existing implementation ran out of memory. Our
implementation was far from having such problems, consuming no more
than about $50$ MB for the Suzuki groups.

The conclusion must therefore be that the project is a small success.

\appendix

\chapter{Benchmark}

Here are details of the results of the first benchmark. The first column shows the
test groups, as they are written in \GAP, and the other columns show
the execution time in milliseconds. The method used to compare the
algorithms was to compute the order of the input group using the orbit
sizes in the computed base and SGS, so to measure the existing
implementation we used a command like
\begin{verbatim}
gap> Size(Group(GeneratorsOfGroup(GL(4, 4))));
gap> benchmark_time := time;
\end{verbatim}

\begin{table}[ht]
\caption{Benchmark results}
\begin{tabular}{l|l|l}
Group & Project & \GAP \\
\hline 
$\SL(2,2)$ & $110$ & $170$ \\
$\GO(+1,2,2)$ & $20$ & $10$ \\
$\GO(-1,2,2)$ & $20$ & $0$ \\
$\GL(2,4)$ & $80$ & $20$ \\
$\SL(2,4)$ & $60$ & $10$ \\
$\GO(+1,2,4)$ &  $40$ & $0$ \\
$\GO(-1,2,4)$ & $40$ & $10$ \\
$\GL(2,3)$ & $80$ & $30$ \\
$\SL(2,3)$ & $50$ & $0$ \\
$\GO(+1,2,3)$ & $50$ & $10$ \\
$\GO(-1,2,3)$ & $50$ & $0$ \\
$\SO(+1,2,3)$ & $20$ & $0$ \\
$\SO(-1,2,3)$ & $30$ & $10$ \\
$\GL(2,5)$ & $70$ & $10$ \\
$\SL(2,5)$ & $40$ & $10$ \\
$\GO(+1,2,5)$ & $50$ & $10$ \\
$\GO(-1,2,5)$ & $50$ & $0$ \\
$\SO(+1,2,5)$ & $50$ & $10$ \\
$\SO(-1,2,5)$ & $40$ & $0$ \\
$\SL(3,2)$ & $50$ & $10$ \\
$\GO(0,3,2)$ & $40$ & $0$ \\
$\GL(3,4)$ & $170$ & $20$ \\
\end{tabular}
\end{table}

\twocolumn
\begin{table}[ht]
\caption{Benchmark results}
\begin{tabular}{l|l|l}
Group & Project & \GAP \\
\hline 
$\SL(3,4)$ & $110$ & $20$ \\
$\GO(0,3,4)$ & $90$ & $10$ \\
$\GL(3,3)$ & $70$ & $10$ \\
$\SL(3,3)$ & $100$ & $0$ \\
$\GO(0,3,3)$ & $60$ & $10$ \\
$\SO(0,3,3)$ & $50$ & $0$ \\
$\GL(3,5)$ & $190$ & $60$ \\
$\SL(3,5)$ & $110$ & $40$ \\
$\GO(0,3,5)$ & $70$ & $10$ \\
$\SO(0,3,5)$ & $70$ & $20$ \\
$\SL(4,2)$ & $130$ & $10$ \\
$\GO(+1,4,2)$ &  $40$ & $10$  \\
$\GO(-1,4,2)$ &  $60$ & $0$ \\
$\GL(4,4)$ & $1480$ & $30$ \\
$\SL(4,4)$ & $690$ & $80$ \\
$\GO(+1,4,4)$ & $90$ & $30$ \\
$\GO(-1,4,4)$ & $170$ & $20$ \\
$\GL(4,3)$ & $210$ & $30$ \\
$\SL(4,3)$ & $410$ & $20$ \\
$\GO(+1,4,3)$ & $80$ & $20$ \\
$\GO(-1,4,3)$ & $80$ & $10$ \\
$\SO(+1,4,3)$ & $80$ & $0$ \\
$\SO(-1,4,3)$ & $80$ & $10$ \\
$\GL(4,5)$ & $1600$ & $50$ \\
$\SL(4,5)$ & $790$ & $190$ \\
$\GO(+1,4,5)$ & $100$ & $60$ \\
$\GO(-1,4,5)$ & $140$ & $50$ \\
$\SO(+1,4,5)$ & $140$ & $60$ \\
$\SO(-1,4,5)$ & $100$ & $60$ \\
$\SL(5,2)$ & $560$ &  $10$ \\
$\GO(0,5,2)$ & $90$ &  $10$ \\
$\GL(5,4)$ & $26440$ & $220$ \\
$\SL(5,4)$ & $4650$ & $440$ \\
$\GO(0,5,4)$ & $230$ & $100$ \\
\end{tabular}
\end{table}

\begin{table}[ht]
\caption{Benchmark results}
\begin{tabular}{l|l|l}
Group & Project & \GAP \\
\hline 
$\GL(5,3)$ & $890$ & $30$ \\
$\SL(5,3)$ & $5090$ & $130$ \\
$\GO(0,5,3)$ & $100$ & $40$ \\
$\SO(0,5,3)$ & $150$ & $40$ \\
$\GL(5,5)$ & $24500$ & $410$ \\
$\SL(5,5)$ & $12380$ & $1450$ \\
$\GO(0,5,5)$ & $470$ & $370$ \\
$\SO(0,5,5)$ & $390$ & $380$ \\
$\SL(6,2)$ & $1250$ & $40$ \\
$\GO(+1,6,2)$ & $250$ & $20$ \\
$\GO(-1,6,2)$ & $350$ & $10$  \\
$\GL(6,4)$ & $205430$ & $460$ \\
$\SL(6,4)$ & $89580$ & $4030$ \\
$\GO(+1,6,4)$ & $3880$ $580$  \\
$\GO(-1,6,4)$ & $3700$ & $550$ \\
$\GL(6,3)$ & $3850$ & $100$ \\
$\SL(6,3)$ & $16890$ & $660$ \\
$\GO(+1,6,3)$ & $340$ & $90$ \\
$\GO(-1,6,3)$ & $340$ & $200$ \\
$\SO(+1,6,3)$ & $520$ & $100$ \\
$\SO(-1,6,3)$ & $480$ & $100$ \\
$\GL(6,5)$ & $292220$ & $84030$ \\
$\SL(6,5)$ & $227830$ & $82010$ \\
$\GO(+1,6,5)$ & $5410$ & $2050$ \\
$\GO(-1,6,5)$ & $1710$ & $1440$ \\
$\SO(+1,6,5)$ & $3190$ & $1460$ \\
$\SO(-1,6,5)$ & $3520$ & $1390$ \\
$\Sz(8)$ & $180$ & $1170$ \\
$\Sz(32)$ & $4390$ & $133860$ \\
$\Sz(128)$ & $1084510$ & $>3600000$ \\
$\Ree(27)$ & $816590$ & $687780$ \\
\end{tabular}
\end{table}
\onecolumn

Here are the results of the second benchmark. The first column shows
the finite field size, the second column the matrix dimension and the
third column is the size of the random generating set that was formed. Each
generating set was computed using the given number of calls to \texttt{RandomInvertibleMat}. For
each line in the table, $20$ generating sets with the given parameters
were formed, the algorithm was executed on these sets, and the average
time over these $20$ executions is shown in the table in the last two
columns.

\begin{table}[ht]
\caption{Benchmark results}
\begin{tabular}{l|l|l|l|l}
Field & Dim & Set & Project & \GAP \\
\hline
$2$ & $2$ & $1$ & $25$ & $16$\\
$2$ & $2$ & $2$ & $23$ & $4$\\
$2$ & $2$ & $3$ & $23$ & $2$\\
$2$ & $2$ & $4$ & $23$ & $2$\\
$2$ & $2$ & $5$ & $23$ & $2$\\
$2$ & $2$ & $6$ & $23$ & $2$\\
$2$ & $2$ & $7$ & $23$ & $3$\\
$2$ & $2$ & $8$ & $23$ & $2$\\
$2$ & $2$ & $9$ & $24$ & $3$\\
$2$ & $2$ & $10$ & $23$ & $3$\\
$2$ & $3$ & $1$ & $34$ & $2$\\
$2$ & $3$ & $2$ & $42$ & $3$\\
$2$ & $3$ & $3$ & $49$ & $3$\\
$2$ & $3$ & $4$ & $52$ & $3$\\
$2$ & $3$ & $5$ & $53$ & $4$\\
$2$ & $3$ & $6$ & $58$ & $4$\\
$2$ & $3$ & $7$ & $60$ & $4$\\
$2$ & $3$ & $8$ & $68$ & $4$\\
$2$ & $3$ & $9$ & $62$ & $5$\\
$2$ & $3$ & $10$ & $68$ & $4$\\
$2$ & $4$ & $1$ & $46$ & $3$\\
$2$ & $4$ & $2$ & $87$ & $4$\\
$2$ & $4$ & $3$ & $117$ & $5$\\
$2$ & $4$ & $4$ & $120$ & $6$\\
$2$ & $4$ & $5$ & $126$ & $6$\\
$2$ & $4$ & $6$ & $124$ & $7$\\
$2$ & $4$ & $7$ & $145$ & $7$\\
$2$ & $4$ & $8$ & $151$ & $7$\\
$2$ & $4$ & $9$ & $161$ & $8$\\
$2$ & $4$ & $10$ & $168$ & $8$\\
\end{tabular}
\end{table}

\twocolumn
\begin{table}[ht]
\caption{Benchmark results}
\begin{tabular}{l|l|l|l|l}
Field & Dim & Set & Project & \GAP \\
\hline
$2$ & $5$ & $1$ & $58$ & $3$\\
$2$ & $5$ & $2$ & $208$ & $10$\\
$2$ & $5$ & $3$ & $231$ & $11$\\
$2$ & $5$ & $4$ & $301$ & $12$\\
$2$ & $5$ & $5$ & $294$ & $13$\\
$2$ & $5$ & $6$ & $357$ & $14$\\
$2$ & $5$ & $7$ & $384$ & $15$\\
$2$ & $5$ & $8$ & $355$ & $15$\\
$2$ & $5$ & $9$ & $407$ & $16$\\
$2$ & $5$ & $10$ & $471$ & $17$\\
$2$ & $6$ & $1$ & $75$ & $4$\\
$2$ & $6$ & $2$ & $612$ & $24$\\
$2$ & $6$ & $3$ & $530$ & $27$\\
$2$ & $6$ & $4$ & $704$ & $28$\\
$2$ & $6$ & $5$ & $1033$ & $32$\\
$2$ & $6$ & $6$ & $903$ & $34$\\
$2$ & $6$ & $7$ & $1209$ & $37$\\
$2$ & $6$ & $8$ & $1055$ & $41$\\
$2$ & $6$ & $9$ & $1300$ & $44$\\
$2$ & $6$ & $10$ & $1158$ & $45$\\
$2$ & $7$ & $1$ & $86$ & $8$\\
$2$ & $7$ & $2$ & $1229$ & $24$\\
$2$ & $7$ & $3$ & $1675$ & $20$\\
$2$ & $7$ & $4$ & $2289$ & $20$\\
$2$ & $7$ & $5$ & $2183$ & $25$\\
$2$ & $7$ & $6$ & $2404$ & $23$\\
$2$ & $7$ & $7$ & $2556$ & $24$\\
$2$ & $7$ & $8$ & $2477$ & $25$\\
$2$ & $7$ & $9$ & $2963$ & $29$\\
$2$ & $7$ & $10$ & $3343$ & $27$\\
$4$ & $2$ & $1$ & $37$ & $3$\\
$4$ & $2$ & $2$ & $55$ & $4$\\
$4$ & $2$ & $3$ & $60$ & $5$\\
$4$ & $2$ & $4$ & $61$ & $5$\\
$4$ & $2$ & $5$ & $62$ & $5$\\
$4$ & $2$ & $6$ & $64$ & $6$\\
$4$ & $2$ & $7$ & $69$ & $7$\\
$4$ & $2$ & $8$ & $70$ & $7$\\
$4$ & $2$ & $9$ & $73$ & $7$\\
$4$ & $2$ & $10$ & $75$ & $8$\\
$4$ & $3$ & $1$ & $59$ & $5$\\
$4$ & $3$ & $2$ & $118$ & $12$\\
$4$ & $3$ & $3$ & $133$ & $14$\\
$4$ & $3$ & $4$ & $152$ & $17$\\
$4$ & $3$ & $5$ & $158$ & $19$\\
$4$ & $3$ & $6$ & $175$ & $21$\\
$4$ & $3$ & $7$ & $185$ & $22$\\
$4$ & $3$ & $8$ & $215$ & $25$\\
$4$ & $3$ & $9$ & $216$ & $27$\\
$4$ & $3$ & $10$ & $226$ & $28$\\
\end{tabular}
\end{table}

\begin{table}[ht]
\caption{Benchmark results}
\begin{tabular}{l|l|l|l|l}
Field & Dim & Set & Project & \GAP \\
\hline
$4$ & $4$ & $1$ & $75$ & $15$\\
$4$ & $4$ & $2$ & $430$ & $26$\\
$4$ & $4$ & $3$ & $576$ & $28$\\
$4$ & $4$ & $4$ & $589$ & $29$\\
$4$ & $4$ & $5$ & $749$ & $28$\\
$4$ & $4$ & $6$ & $714$ & $31$\\
$4$ & $4$ & $7$ & $838$ & $33$\\
$4$ & $4$ & $8$ & $912$ & $38$\\
$4$ & $4$ & $9$ & $1019$ & $36$\\
$4$ & $4$ & $10$ & $1085$ & $38$\\
$4$ & $5$ & $1$ & $109$ & $52$\\
$4$ & $5$ & $2$ & $3141$ & $170$\\
$4$ & $5$ & $3$ & $4248$ & $119$\\
$4$ & $5$ & $4$ & $3660$ & $106$\\
$4$ & $5$ & $5$ & $5126$ & $116$\\
$4$ & $5$ & $6$ & $3895$ & $124$\\
$4$ & $5$ & $7$ & $4416$ & $134$\\
$4$ & $5$ & $8$ & $5154$ & $137$\\
$4$ & $5$ & $9$ & $5002$ & $147$\\
$4$ & $5$ & $10$ & $6285$ & $156$\\
$4$ & $6$ & $1$ & $193$ & $212$\\
$4$ & $6$ & $2$ & $27943$ & $655$\\
$4$ & $6$ & $3$ & $33664$ & $455$\\
$4$ & $6$ & $4$ & $30882$ & $478$\\
$4$ & $6$ & $5$ & $40467$ & $506$\\
$4$ & $6$ & $6$ & $36702$ & $534$\\
$4$ & $6$ & $7$ & $35179$ & $570$\\
$4$ & $6$ & $8$ & $41211$ & $593$\\
$4$ & $6$ & $9$ & $47733$ & $630$\\
$4$ & $6$ & $10$ & $46698$ & $661$\\
$4$ & $7$ & $1$ & $8786$ & $19238$\\
$4$ & $7$ & $2$ & $335656$ & $95870$\\
$4$ & $7$ & $3$ & $337232$ & $100999$\\
$4$ & $7$ & $4$ & $371067$ & $98525$\\
$4$ & $7$ & $5$ & $381631$ & $94616$\\
$4$ & $7$ & $6$ & $463195$ & $94353$\\
$4$ & $7$ & $7$ & $463824$ & $102045$\\
$4$ & $7$ & $8$ & $460079$ & $104262$\\
$4$ & $7$ & $9$ & $485832$ & $115177$\\
$4$ & $7$ & $10$ & $497838$ & $99240$\\
$3$ & $2$ & $1$ & $39$ & $3$\\
$3$ & $2$ & $2$ & $47$ & $3$\\
$3$ & $2$ & $3$ & $49$ & $3$\\
$3$ & $2$ & $4$ & $52$ & $3$\\
$3$ & $2$ & $5$ & $51$ & $4$\\
$3$ & $2$ & $6$ & $53$ & $4$\\
$3$ & $2$ & $7$ & $46$ & $4$\\
$3$ & $2$ & $8$ & $54$ & $4$\\
$3$ & $2$ & $9$ & $54$ & $5$\\
$3$ & $2$ & $10$ & $47$ & $5$\\
\end{tabular}
\end{table}

\begin{table}[ht]
\caption{Benchmark results}
\begin{tabular}{l|l|l|l|l}
Field & Dim & Set & Project & \GAP \\
\hline
$3$ & $3$ & $1$ & $59$ & $4$\\
$3$ & $3$ & $2$ & $92$ & $7$\\
$3$ & $3$ & $3$ & $98$ & $8$\\
$3$ & $3$ & $4$ & $101$ & $9$\\
$3$ & $3$ & $5$ & $111$ & $9$\\
$3$ & $3$ & $6$ & $114$ & $11$\\
$3$ & $3$ & $7$ & $126$ & $12$\\
$3$ & $3$ & $8$ & $126$ & $13$\\
$3$ & $3$ & $9$ & $136$ & $12$\\
$3$ & $3$ & $10$ & $143$ & $15$\\
$3$ & $4$ & $1$ & $74$ & $7$\\
$3$ & $4$ & $2$ & $244$ & $22$\\
$3$ & $4$ & $3$ & $278$ & $26$\\
$3$ & $4$ & $4$ & $310$ & $29$\\
$3$ & $4$ & $5$ & $366$ & $32$\\
$3$ & $4$ & $6$ & $334$ & $35$\\
$3$ & $4$ & $7$ & $415$ & $39$\\
$3$ & $4$ & $8$ & $452$ & $40$\\
$3$ & $4$ & $9$ & $475$ & $45$\\
$3$ & $4$ & $10$ & $454$ & $47$\\
$3$ & $5$ & $1$ & $82$ & $15$\\
$3$ & $5$ & $2$ & $1029$ & $44$\\
$3$ & $5$ & $3$ & $1110$ & $28$\\
$3$ & $5$ & $4$ & $1692$ & $29$\\
$3$ & $5$ & $5$ & $1361$ & $36$\\
$3$ & $5$ & $6$ & $1612$ & $34$\\
$3$ & $5$ & $7$ & $1752$ & $35$\\
$3$ & $5$ & $8$ & $1887$ & $40$\\
$3$ & $5$ & $9$ & $1928$ & $39$\\
$3$ & $5$ & $10$ & $2197$ & $42$\\
$3$ & $6$ & $1$ & $111$ & $44$\\
$3$ & $6$ & $2$ & $6517$ & $105$\\
$3$ & $6$ & $3$ & $5772$ & $131$\\
$3$ & $6$ & $4$ & $7463$ & $117$\\
$3$ & $6$ & $5$ & $6184$ & $104$\\
$3$ & $6$ & $6$ & $7111$ & $128$\\
$3$ & $6$ & $7$ & $7191$ & $117$\\
$3$ & $6$ & $8$ & $8478$ & $120$\\
$3$ & $6$ & $9$ & $11749$ & $127$\\
$3$ & $6$ & $10$ & $10142$ & $134$\\
$3$ & $7$ & $1$ & $385$ & $128$\\
$3$ & $7$ & $2$ & $33440$ & $591$\\
$3$ & $7$ & $3$ & $36465$ & $586$\\
$3$ & $7$ & $4$ & $35244$ & $506$\\
$3$ & $7$ & $5$ & $35492$ & $331$\\
$3$ & $7$ & $6$ & $43913$ & $351$\\
$3$ & $7$ & $7$ & $47698$ & $364$\\
$3$ & $7$ & $8$ & $45104$ & $382$\\
$3$ & $7$ & $9$ & $46772$ & $401$\\
$3$ & $7$ & $10$ & $62226$ & $421$\\
\end{tabular}
\end{table}

\begin{table}[ht]
\caption{Benchmark results}
\begin{tabular}{l|l|l|l|l}
Field & Dim & Set & Project & \GAP \\
\hline
$5$ & $2$ & $1$ & $46$ & $3$\\
$5$ & $2$ & $2$ & $59$ & $5$\\
$5$ & $2$ & $3$ & $61$ & $6$\\
$5$ & $2$ & $4$ & $64$ & $6$\\
$5$ & $2$ & $5$ & $64$ & $8$\\
$5$ & $2$ & $6$ & $65$ & $9$\\
$5$ & $2$ & $7$ & $68$ & $9$\\
$5$ & $2$ & $8$ & $73$ & $10$\\
$5$ & $2$ & $9$ & $72$ & $10$\\
$5$ & $2$ & $10$ & $77$ & $11$\\
$5$ & $3$ & $1$ & $59$ & $8$\\
$5$ & $3$ & $2$ & $147$ & $14$\\
$5$ & $3$ & $3$ & $157$ & $15$\\
$5$ & $3$ & $4$ & $178$ & $15$\\
$5$ & $3$ & $5$ & $182$ & $14$\\
$5$ & $3$ & $6$ & $203$ & $15$\\
$5$ & $3$ & $7$ & $226$ & $17$\\
$5$ & $3$ & $8$ & $240$ & $17$\\
$5$ & $3$ & $9$ & $258$ & $19$\\
$5$ & $3$ & $10$ & $269$ & $20$\\
$5$ & $4$ & $1$ & $78$ & $35$\\
$5$ & $4$ & $2$ & $957$ & $94$\\
$5$ & $4$ & $3$ & $1123$ & $73$\\
$5$ & $4$ & $4$ & $1038$ & $65$\\
$5$ & $4$ & $5$ & $1300$ & $68$\\
$5$ & $4$ & $6$ & $1337$ & $72$\\
$5$ & $4$ & $7$ & $1563$ & $78$\\
$5$ & $4$ & $8$ & $1703$ & $81$\\
$5$ & $4$ & $9$ & $1815$ & $93$\\
$5$ & $4$ & $10$ & $1864$ & $91$\\
$5$ & $5$ & $1$ & $166$ & $158$\\
$5$ & $5$ & $2$ & $10330$ & $414$\\
$5$ & $5$ & $3$ & $11236$ & $508$\\
$5$ & $5$ & $4$ & $11779$ & $402$\\
$5$ & $5$ & $5$ & $13499$ & $355$\\
$5$ & $5$ & $6$ & $15476$ & $375$\\
$5$ & $5$ & $7$ & $14578$ & $475$\\
$5$ & $5$ & $8$ & $14204$ & $428$\\
$5$ & $5$ & $9$ & $14994$ & $450$\\
$5$ & $5$ & $10$ & $15487$ & $472$\\
$5$ & $6$ & $1$ & $903$ & $3568$\\
$5$ & $6$ & $2$ & $180644$ & $86753$\\
$5$ & $6$ & $3$ & $195962$ & $91533$\\
$5$ & $6$ & $4$ & $210412$ & $84445$\\
$5$ & $6$ & $5$ & $216448$ & $84051$\\
$5$ & $6$ & $6$ & $225522$ & $80495$\\
$5$ & $6$ & $7$ & $212674$ & $76792$\\
$5$ & $6$ & $8$ & $258668$ & $81668$\\
$5$ & $6$ & $9$ & $269393$ & $97520$\\
$5$ & $6$ & $10$ & $287699$ & $87331$
\end{tabular}
\end{table}

\onecolumn

\backmatter

\bibliographystyle{amsalpha}
\bibliography{schreiersims}

\end{document}